\documentclass[journal,onecolumn]{IEEEtran} %STIX1COL,STIX2COL,STIXSMALL

\usepackage{algorithmic}
\usepackage{graphicx}
\usepackage{textcomp}
\def\BibTeX{{\rm B\kern-.05em{\sc i\kern-.025em b}\kern-.08em
    T\kern-.1667em\lower.7ex\hbox{E}\kern-.125emX}}

\usepackage{umlaute,eurosym,comment}
\usepackage{amsmath,amssymb,amsfonts,bm,upgreek}
\usepackage{float,epsfig,pstricks,pst-node,subfig}
\usepackage{wrapfig}
\usepackage{ifpdf,cite}
\usepackage{multirow,dcolumn}
\usepackage{tikz,pgfplots}
\usepackage{subcaption}

\psset{unit=1mm,framesep=10pt,fillstyle=none}  % pstricks basic settings
\degrees[360]                                  % segmentation of the unit circle
\psset{linewidth=0.9pt}                        % standard line width
\graphicspath{{figures/}}

\usetikzlibrary{arrows,chains,matrix,positioning,scopes}                               
\usetikzlibrary{intersections, pgfplots.fillbetween}

\newlength\figureheight	
\newlength\figurewidth	

\newtheorem{theorem}{Theorem}
\newtheorem{lemma}[theorem]{Lemma}

\newtheorem{corollary}[theorem]{Corollary}

\newtheorem{definition}[theorem]{Definition}
\newtheorem{remark}[theorem]{Remark}
\newtheorem{assumption}[theorem]{Assumption}
\newtheorem{ROBass}[theorem]{Assumption R.\!}
\newtheorem{LPVass}[theorem]{Assumption L.\!}
\providecommand{\ef}{\;.}
\providecommand{\ec}{\;,}

\providecommand{\fa}{\forall\;}         % for all
\providecommand{\eqdef}{\triangleq}     % equals by definition
\providecommand{\eqcon}{\overset{!}{=}} % equals by condition
\providecommand{\to}{\rightarrow}                        % converges to
\newcommand{\conv}[1]{\text{conv}\left\{#1\right\}}        % convex hull
\newcommand{\vct}[1]{\mathrm{vec}\left(#1\right)}  % vectorization

\newcommand{\BZ}{\mathbb{Z}}        % integers
\newcommand{\BR}{\mathbb{R}}        % real numbers
\newcommand{\BsX}{\mathbb{X}} 
\newcommand{\BU}{\mathbb{U}}
\newcommand{\BD}{\mathbb{D}} 
\newcommand{\BDt}{\mathbb{D}_\theta} 
\newcommand{\CX}{\mathcal{X}} 
\newcommand{\CP}{\mathcal{P}} 
\newcommand{\BW}{\mathbb{W}} 

\definecolor{ULOcean}{RGB}{0,75,90}
\definecolor{ULOrange}{RGB}{236,116,4}
\definecolor{LiteRed}{RGB}{255,206,206}
\definecolor{MedRed}{RGB}{204,24,24}
\definecolor{DarkRed}{RGB}{127,15,15}
\definecolor{LiteYellow}{RGB}{255,255,100}
\definecolor{MedYellow}{RGB}{255,255,25}
\definecolor{DarkYellow}{RGB}{225,225,0}
\definecolor{LiteBlue}{RGB}{170,184,255}
\definecolor{MedBlue}{RGB}{10,20,204}
\definecolor{DarkBlue}{RGB}{27,29,120}
\definecolor{LiteGreen}{RGB}{24,204,75}
\definecolor{MedGreen}{RGB}{0,130,37}
\definecolor{DarkGreen}{RGB}{0,61,17}
\definecolor{LiteViolet}{RGB}{110,0,225}
\definecolor{MedViolet}{RGB}{73,0,148}
\definecolor{DarkViolet}{RGB}{53,0,107}
\definecolor{LiteOrange}{RGB}{255,226,117}
\definecolor{MedOrange}{RGB}{255,206,0}
\definecolor{DarkOrange}{RGB}{204,167,10}
\definecolor{LiteGray}{RGB}{230,230,230}
\definecolor{MedGray}{RGB}{160,160,160}
\definecolor{DarkGray}{RGB}{100,100,100}

\begin{document}
\title{Deadbeat Robust Model Predictive Control for Linear Parameter-Varying Systems}
\author{Georg Schildbach, \IEEEmembership{Member, IEEE}
\thanks{Submitted to International Journal of Robust and Nonlinear Control on September 16, 2025. This work has been supported by the German Research Foundation (DFG) under grants no.\ 460891204 and no.\ 419290163.}
\thanks{Georg Schildbach is with the University of L{\"u}beck, Ratzeburger Allee 160, 23568 L{\"u}beck, Germany (e-mail: georg.schildbach@uni-luebeck.de).}
Hossam S.\ Abbas, \IEEEmembership{Member, IEEE}
\thanks{Hossam S.\ Abbas is with the University of L{\"u}beck, Ratzeburger Allee 160, 23568 L{\"u}beck, Germany (e-mail: h.abbas@uni-luebeck.de).}}

\maketitle

\begin{abstract}
The concept of Deadbeat Robust Model Predictive Control (DRMPC)  is to completely extinguish the effect of external disturbances within the first few steps of the prediction horizon. The benefit is that the remaining dynamics of the system can be planned with certainty. This means it is not necessary to employ a Robust Positive Invariant (RPI) set as a terminal condition, which is often complex or even intractable to compute. Recent work has shown that it is possible to obtain full system theoretic guarantees in the case of linear systems, including recursive feasibility of all constraints and input-to-state stability. This paper extends this contribution to linear time-varying (LTV) or parameter-varying systems (LPV) systems with bounded parametric uncertainty and additive disturbances. Full system-theoretic guarantees can also be provided for these cases. Numerical simulation results demonstrate that the performance of the proposed LPV-DRMPC scheme, with the origin as the terminal set and a slightly increased prediction horizon, is almost comparable to that of other LPV-MPC schemes with an RPI terminal set constraint, despite a lower computational complexity. Besides avoiding the computation of RPI terminal sets, LPV-DRMPC allows to shift a significant portion of the computations offline.
\end{abstract}

\begin{IEEEkeywords}
Model Predictive Control, linear time-varying systems, linear parameter-varying systems, robust control, invariant sets
\end{IEEEkeywords}

%\jnlcitation{\cname{%
%\author{Schildbach G.}.
%\ctitle{On simplifying ‘incremental remap’-based transport schemes.} \cjournal{\it J Comput Phys.} \cvol{2021;00(00):1--18}.}

\maketitle

\vspace*{-0.6cm}
\section{Introduction}\label{Sec:Intro}

Linear Parameter-Varying Model Predictive Control (LPV-MPC) is a recent and powerful control strategy for nonlinear, time-varying systems, and systems  subject to time-varying  parametric uncertainty \cite{Morato:2020}. It aims to account for the parametric uncertainty and essential nonlinearity of the plant via their embedding into a Linear Parameter-Varying (LPV) system. The goal is to reduce the computational requirements compared to a full-blown nonlinear approach, by applying an extension of the classic methodology of linear Model Predictive Control (MPC). LPV-MPC has demonstrated an excellent performance in a wide range of practical applications, including automated vehicles \cite{AlcalaEtAl:2020,NezamiLPV:2022,KaraAbbas:2025}, aerospace systems \cite{CavaniniEtAl:2021,JinEtAl:2021,BrysonGruenwald:2022}, and wind turbines \cite{HanGao:2022,LiuEtAl:2025}.

From a research perspective, several challenges related to LPV-MPC remain open. Some of them are related to specific applications, such as the identification and proper formulation of LPV models\cite{Toth:2010} and others are realted to the design of proper real-time implementations using explicit solutions or online optimization solvers \cite{Hespe:2021, KaraAbbas:2025}. Also the fundamental design of LPV-MPC algorithms is subject to ongoing research \cite{Abbas:2019, Hanema:2021,Heydari:2021,Abbas:2024}. A major challenge is to find approaches that are robust to parameter-related uncertainties over the MPC prediction horizon, without being overly conservative or highly computationally demanding due to the online construction of robustification conditions, e.g., invariant tubes, constraint tightening, or computation of robust positive invariant sets. These approaches require rigorous proofs of recursive feasibility and robust stability. Moreover, they should be easy to implement and computationally efficient. %\textcolor{red}{Hossam: Maybe you can add a few more references into this paragraph?}

The contribution of this paper is related to the fundamental design of LPV-MPC algorithms. It leverages the recent results on Deadbeat Robust Model Predictive Control (DRMPC) for linear systems \cite{Schildbach:2025}. The main idea of DRMPC is to deadbeat the effect of the uncertainty over the first few steps of the prediction horizon, which is called the \emph{deadbeat horizon}. This deadbeat property inside the predictions is only virtual in nature, in contrast to actual deadbeat control \cite{EmamFran:1982}. The advantage lies in the complete elimination of the uncertainty in the predictions beyond the deadbeat horizon. This allows for the use of a deterministic terminal condition, which significantly simplifies the controller design and implementation, while providing the full range of system-theoretic guarantees. The control performance of DRMPC appears competitive with other state-of-the-art approach of Robust Model Predictive Control (RMPC) for linear systems \cite{Schildbach:2025}, such as Tube-based RMPC \cite{Mayne:2005,RacEtAl:2012,RaMaDi:2018} and RMPC with Affine Disturbance Feedback \cite{Bemporad:1998,Loefberg:2003a,Goulart:2006}.

The contributions of this paper can be summarized as follows: (1) introduction of a new parameterization of the deadbeat predictions for linear and LPV systems; (2) extension of DRMPC to a class of LPV systems; (3) rigorous proofs of recursive feasibility and robust exponential stability for the linear and the LPV Case; (4) evaluation of the performance based on a commonly used benchmark example \cite{FlemingEtAl:2015}.

\section{Background}

\subsection{Terminology and Notation}

$\mathbb{R}$ and $\mathbb{Z}$ denote the sets of real and integral numbers. $\mathbb{R}_{+}$ ($\mathbb{R}_{0+}$) and $\mathbb{Z}_{+}$ ($\mathbb{Z}_{0+}$) are then the sets of positive (non-negative) real and integral numbers, respectively. 

$\mathbb{R}^{n}$ represents a vector space of dimension $n\in\mathbb{Z}_{+}$. The $i$-th component of a vector $v\in\mathbb{R}^{n}$ is denoted $v[i]\in\BR$. The notation $\|\cdot\|$ is used to denote any norm. If $A\in\BR^{n\times n}$ is a matrix, $\vct{A}\in\BR^{n^2}$ denotes the vectorization of $A$, i.e., a column vector containing the stacked-up columns of $A$.
The notation $\conv{\cdot}$ denotes the convex hull of a set.

If $\mathbb{S},\mathbb{T}\subset\mathbb{R}^{n}$ are sets, the \textbf{Minkowski sum} is defined as
\begin{equation*}
	\mathbb{S}\oplus \mathbb{T}\triangleq\bigl\{s+t\in\mathbb{R}^{n}\:\big|\: s\in \mathbb{S},\:t\in \mathbb{T}\}\ec
\end{equation*}
and the \textbf{Pontryagin difference} as
\begin{equation*}
	\mathbb{S}\ominus \mathbb{T}\triangleq\bigl\{\xi\in\mathbb{R}^{n}\:\big|\:\xi+t\in \mathbb{S}\:\:\forall\:t\in \mathbb{T}\}\;.
\end{equation*}

A polytope is a compact polyhedron, which is the inter-
section of a finite number of half spaces. A (hyper)box
is a  polytope where all the defining hyperplanes
are axis parallel.

A function $\alpha:\BR_{0+}\rightarrow\BR_{0+}$ is a \textbf{$\textbf{K}$-function} if it is continuous, strictly monotonically increasing, and $\alpha(0)=0$. It is a \textbf{$\textbf{K}_{\infty}$-function} if, in addition, $\alpha(r)\rightarrow\infty$ as $r\rightarrow\infty$. A function $\beta:\BR_{0+}\times\BZ_{0+}\rightarrow\BR_{0+}$ is a \textbf{$\textbf{KL}$-function} if $\beta(\,\cdot\,,k)$ is a $\text{K}$-function for any fixed $k\in\BZ_{0+}$, and $\beta(r,\,\cdot\,)$ is monotonically decreasing with $\beta(r,k)\rightarrow 0$ as $k\rightarrow\infty$ for any fixed $r\in\BR_{0+}$.

\subsection{Control System}

This paper considers a discrete-time (DT) linear state space model
\begin{equation}\label{Equ:DTSystem}
  x_{k+1}=A(\theta_{k})x_{k}+B(\theta_{k})u_{k}+w_{k}
\end{equation}
that depends on a time-varying parameter $\theta_{k}\in\mathbb{R}^{p}$. Here $k\in\mathbb{Z}_{0+}$ denotes the time step, $x_{k}\in\mathbb{R}^{n}$ is the state of the system, $u_{k}\in\mathbb{R}^{m}$ is the control input, and $w_{k}\in\mathbb{R}^{n}$ is an additive disturbance. The system comes with a corresponding initial condition $x_{0}\in\mathbb{R}^{n}$. The additive disturbance and the parametric uncertainty are time-varying, yet contained in constant and known sets
\begin{subequations}\label{Equ:ParSet}\begin{align}
  w_{k}&\in\mathbb{W}\subseteq\BR^{n}\quad\fa k\in\mathbb{Z}_{0+}\ec\\
  \theta_{k}&\in\Theta\subseteq\BR^{p}\quad\fa k\in\mathbb{Z}_{0+}\;
\end{align}\end{subequations}
where $\mathbb{W}, \Theta$ are compact sets.

\begin{assumption}[measurements of the state]\label{The:Measurement-x}
  The state of the system $x_{k}$ is measured in each step $k$. 
\end{assumption}

This paper considers two cases.
In the first case, we assume that $\theta_k$ is completely arbitrary (within the compact set $\Theta$). This is referred to as the \emph{Robust Case}. The system matrix $A(\theta_{k}) \in \mathbb{R}^{n \times n}$ and the input matrix $B(\theta_{k}) \in \mathbb{R}^{n \times m}$ are assumed to depend affinely on the time-varying parameter,
\begin{subequations}\begin{align}\label{Equ:Parameter}
  A(\theta_{k})&=\bar{A}+\sum_{i=1}^{p}\theta_{k}[i]A^{(i)}\ec\\
  B(\theta_{k})&=\bar{B}+\sum_{i=1}^{p}\theta_{k}[i]B^{(i)}\ef
\end{align}\end{subequations}
Here $\bar{A}$, $A^{(i)}\in\BR^{n\times n}$ and $\bar{B}$, $B^{(i)}\in\BR^{n\times m}$ are constant matrices, where $i\in\{1,2,\dots,p\}$.
\begin{ROBass}[controllability in the Robust Case]\label{The:Controllability}
  The nominal system matrix $\bar{A}$ and the nominal input matrix $\bar{B}$ form a controllable pair.
\end{ROBass}

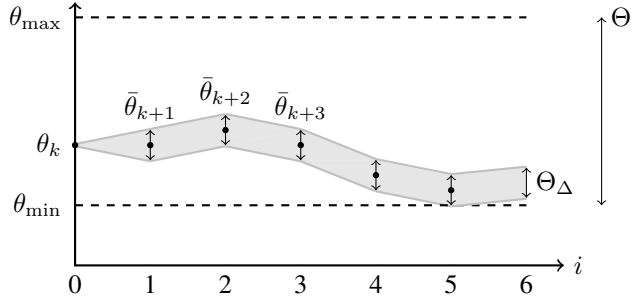
\begin{figure}
\centering
\scalebox{1}{\begin{tikzpicture}[scale=1]
\colorlet{ColorGray}{gray!40}
    \draw [<->,thick] (0,3.5) node (yaxis) [above] {}
        |- (6.5,0) node (xaxis) [right] {$i$};
    \node [left] at (-0.05,1.6) {$\theta_k$}; 
    % Draw two intersecting lines
    \draw  [gray!50,  ultra  thick, name path=c] (0,1.6) -- (2,2.0);
   \draw  [ gray!50, ultra thick ,name path=d] (0,1.6) coordinate (b_1) -- (1,1.4) coordinate (b_2);
    \draw  [ gray!50, ultra thick, name path=e] (2,2) -- (3,1.8);
   \draw  [ gray!50, ultra thick ,name path=f] (1,1.4)  -- (2,1.6);
    \draw  [gray!50, ultra thick, name path=g] (3,1.8) -- (4,1.4);
   \draw  [ gray!50, ultra thick ,name path=h] (2,1.6)  -- (3,1.4);
   \draw  [gray!50, ultra  thick ,name path=i] (3,1.4)  -- (4,1);
   \draw  [gray!50, ultra  thick ,name path=j] (4,1.4)  -- (5,1.2);
   \draw  [ gray!50, ultra thick ,name path=k] (4,1)  -- (5,0.8);
\tikzfillbetween[of=c and d ]{gray!20};
\tikzfillbetween[of=e and f ]{gray!20};
\tikzfillbetween[of=g and h ]{gray!20};
\tikzfillbetween[of=i and j ]{gray!20};
\tikzfillbetween[of=j and k ]{gray!20};
   \draw [<->]   (1,1.8) -- (1,1.4);
   \draw [<->]   (2,2)   -- (2,1.6);
          \draw [<->]   (3,1.8) -- (3,1.4);
          \draw [<->]   (4,1.4) -- (4,1);
\draw [dashed, thick] (0,3.3) -- (6,3.3);  
    \node [left] at (-0.05,3.3) {$\theta_{\max}$}; 
\draw [dashed, thick] (0,0.8) -- (6,0.8);
       \node [left] at (-0.05,0.8) {$\theta_{\min}$};    
   \foreach \x in {0,...,6}
\draw (\x, 0) -- (\x, 0) node [below] {\x};
    %\draw  [thick, name path=c] (1,1.6) -- (2,1.8);
%   \draw  [dashed, red, thick ,name path=d] (1,1.65) coordinate (b_1) -- (2,1.6) coordinate (b_2);
%    \draw  [dashed, red,thick, name path=e] (2,2) -- (3,1.8);
  % \draw  [dashed, red,thick ,name path=f] (1,1.4)  -- (2,1.6);
%    \draw  [dashed, red,thick, name path=g] (3,1.8) -- (4,1.4);
%   \draw  [dashed, red,thick ,name path=h] (2,1.6)  -- (3,1.4);
%   \draw  [dashed, red,thick ,name path=i] (3,1.4)  -- (4,1);
%   \draw  [dashed, red,thick ,name path=j] (4,1.4)  -- (5,1.2);
%   \draw  [dashed, red,thick ,name path=k] (4,1)  -- (5,0.8);
   \draw  [gray!50, ultra  thick , ,name path=j] (5,1.2)  -- (6,1.3);
   \draw  [gray!50, ultra  thick , ,name path=k] (5,0.8)  -- (6,0.9);
\tikzfillbetween[of=j and k ]{gray!20};
   \draw [<->]   (5,1.2) -- (5,0.8);
      \draw [<->]   (6,1.3) -- (6,0.9);  
          \draw [<->]   (7,0.8) -- (7,3.3); 
    \node [right] at (7,3.3) {$\Theta$}; 
\draw [fill] (0,1.6) circle [radius=0.035];
    \draw [fill] (1,1.6) circle [radius=0.035];
        \draw [fill] (2,1.8) circle [radius=0.035];
        \draw [fill] (3,1.6) circle [radius=0.035];
                \draw [fill] (4,1.2) circle [radius=0.035];
                \draw [fill] (5,1.) circle [radius=0.035];
                %\draw [fill] (6,1.1) circle [radius=0.035];
       \node [right] at (6,1.1) {$\Theta_\Delta$};  
      \node [above] at (1,1.8) {$\bar{\theta}_{k+1}$};  
      \node [above] at (2,2.0) {$\bar{\theta}_{k+2}$};  
      \node [above] at (3,1.8) {$\bar{\theta}_{k+3}$};  
\end{tikzpicture}} 
\caption{Illustration of Assumption L.\ref{The:Measurement-theta}. Example of a nominal sequence $\bar{\theta}_{k+j}$ (solid black circles) of a scalar parameter $\theta_{k}$, starting at step $k$ over six prediction steps $j=1,2,\dots,6$. Here $\bar{\theta}_k = \theta_k$ and $\theta_{k+j} = \bar{\theta}_{k+j} + \Delta\theta_{k+j}$ with $\Delta\theta_{k+j} \in \Theta_\Delta$ (gray area), and $\bar{\theta}_{k+j} + \Delta\theta_{k+j}\subseteq\Theta$ (dashed lines) for all $j=1,2,\dots,6$.}
\label{Fig:schtube}
\end{figure}

In the second case, additionally, we assume that $\theta_k$ can be predicted based on some known parameter dynamics, subject to a bounded uncertainty. This is referred to as the \emph{LPV Case}.
\begin{LPVass}[parameter measurements in the LPV Case]\label{The:Measurement-theta}
(a) The time-varying parameter $\theta_{k}$ is measured at each step $k$. (b) Based on $\theta_{k}$ in step $k$, at least $N$ future values of the parameter can be predicted with a (given) bounded uncertainty, denoted $\Theta_\Delta\subset\mathbb{R}^p$. (c) The uncertainties are centered around the predicted nominal values $\bar{\theta}_{k+j}\in\mathbb{R}^p$, $j=1,\dots,N$. 
\end{LPVass}
\noindent Assumption L.\ref{The:Measurement-theta} is illustrated in Fig.~\ref{Fig:schtube}. As a consequence, the matrices $A(\theta_{k})$ and $B(\theta_{k})$ are completely known in step $k$. Their future values, however, are unknown, and assumed to depend affinely on the prediction uncertainties $\Delta\theta_{k+j}\in\Theta_\Delta$ as follows:
\begin{subequations}\begin{align}\label{Equ:Parameter-LPV}
  A(\theta_{k+j})&=\bar{A}_{k+j}+\sum_{i=1}^{p}\Delta\theta_{k+j}[i]A^{(i)}\ec\\
  B(\theta_{k+j})&=\bar{B}_{k+j}+\sum_{i=1}^{p}\Delta\theta_{k+j}[i]B^{(i)}\ef
\end{align}\end{subequations}
Here $\bar{A}_{k+j}\eqdef A(\bar{\theta}_{k+j})$ and $\bar{B}_{k+j}\eqdef B(\bar{\theta}_{k+j})$ are known parameter-dependent matrices for all $j\in\{1,2,\dots,N\}$, and $A^{(i)}, B^{(i)}$ are known matrices for all $i\in\{1,2,\dots,p\}$.
  %\begin{assumption}[measurements of the state]\label{The:Measurement-theta}The uncertain parameter $\theta_{k}$ is measured in each step $k$, and its future values, for at least $N$ steps ahead, is known within a given bounded error $\Delta\theta$; see Fig.~\ref{Fig:schtube} for an illustration.
%\end{assumption}

\begin{LPVass}[controllability in the LPV Case]\label{The:ControllabilityLPV}
 For the nominal system in the LPV Case,
 \begin{equation}\label{Equ:NominalSystemLPV}
     x_{k+1}=\bar{A}_k x_k+\bar{B}_k u_k\ec
 \end{equation}
the matrix
 \begin{equation}\label{Equ:ControllabilityMatrixLPV-1}
     \mathcal{C}_{k,r} = 
\begin{bmatrix}
\bar{B}_{k+r-1} & \bar{A}_{k+r-1} \bar{B}_{k+r-2} & \bar{A}_{k+r-1}\bar{A}_{k+r-2} \bar{B}_{k+r-3} & \cdots & \left(\prod_{j=1}^{r-1} \bar{A}_{k+j}\right) \bar{B}_{k}
\end{bmatrix}
 \end{equation}
  has full rank $n$,  where 
 \[
\prod_{j=1}^{r-1} \bar{A}_{k+j} \triangleq \bar{A}_{k+r-1} \bar{A}_{k+r-2} \cdots \bar{A}_{k+1}\ec
\]
with $r\in\mathbb{Z}_{+}$ such that $r\leq n$.
\end{LPVass}
\noindent
Assumption L.\ref{The:ControllabilityLPV} establishes a controllability condition for the nominal LPV system~\eqref{Equ:NominalSystemLPV} over the interval $[k, k + r]$, where $\mathcal{C}_{k,r}$ is the \emph{controllability matrix} \cite{Kamen:2010}.

%The system matrix $A(\theta_{k})\in\BR^{n\times n}$ and the input matrix $B(\theta_{k})\in\BR^{n\times m}$ affinely depend on the time-varying parameter,
%\begin{subequations}\begin{align}\label{Equ:Parameter}
%  A(\theta_{k})&=\bar{A}+\sum_{i=1}^{p}\theta_{k}[i]A^{(i)}\ec\\
%  B(\theta_{k})&=\bar{B}+\sum_{i=1}^{p}\theta_{k}[i]B^{(i)}\ef
%\end{align}\end{subequations}
%Here $\bar{A}$, $A^{(i)}\in\BR^{n\times n}$ and $\bar{B}$, $B^{(i)}\in\BR^{n\times m}$ are constant matrices, where $i\in\{1,2,\dots,p\}$.

The control objective in both cases is to regulate the state of system \eqref{Equ:DTSystem} to the origin, as good as possible, regardless of the uncertainties and disturbances. As part of the control task, state and input constraints must be respected:
\begin{equation}\label{Equ:Constraints}
  x_{k}\in\BsX\subseteq\BR^{n}\ec\quad u_{k}\in\BU\subseteq\BR^{m}\qquad\fa k\in\mathbb{Z}_{0+}\ef
\end{equation}
\begin{assumption}[constraint and uncertainty sets]\label{The:Constraints}
  (a) The input and state constraint sets $\BU$ and $\BsX$ are convex and contain the origin. (b) The parameter set $\Theta$ is a convex polytope that contains the origin. (c) The disturbance set $\mathbb{W}$ is a convex polytope that contains the origin. (d) The  set $\Theta_\Delta$ is a convex polytope that contains the origin.
\end{assumption}
%%---------------------------------------------------------------
\subsection{Conversion to Additive Uncertainty}\label{Sec:AdditiveUncertainty}

%[Describing the conversion from multiplicative into additive uncertainty]

\subsubsection*{The Robust Case}

The goal of this section is to convert system \eqref{Equ:DTSystem} into the form
\begin{equation}\label{Equ:SystemWAdditiveDis}
x_{k+1}=\bar{A}x_{k}+\bar{B}u_{k}+d_{k}\ec
\end{equation}
where $d_{k}\in\BD$ is a modified disturbance that subsumes the original additive disturbance and the parametric uncertainty. The representation \eqref{Equ:SystemWAdditiveDis} is achieved by defining
\[
d_k\triangleq \left(\sum_{i=1}^{p}\theta_{k}[i]A^{(i)}\right)x_k+\left(\sum_{i=1}^{p}\theta_{k}[i]A^{(i)}\right)u_k+w_k\ef
\]
The set $\mathbb{D}$ is defined by
\begin{equation}
    \BD\triangleq \conv{\left.\left(\sum_{i=1}^{p}\theta_{k}[i]A^{(i)}\right)x_k+\left(\sum_{i=1}^{p}\theta_{k}[i]A^{(i)}\right)u_k +w_k \right| \theta_k\in\Theta\ec w_k\in\BW\ec x_k\in\BsX\ec u_k\in\BU}\ef
\end{equation}
Clearly, under Assumption~\ref{The:Constraints}, $\BD$ is again a convex polytope that contains the origin. Hence the problem statement fully corresponds to the original formulation of DRMPC \cite{Schildbach:2025}.

\subsubsection*{The LPV Case}

At step $k$, only the additive disturbance $w_k \in \mathbb{W}$ acts on the system, since the time-varying parameter $\theta_k$ is known according to Assumption~L.\ref{The:Measurement-theta}. Therefore, the system is represented by
\begin{equation}\label{Equ:LPVSystemWAdditiveDisOnly}
x_{k+1} = A_k x_k + B_k u_k + w_k\ec
\end{equation}
where $A_k = A(\theta_k)$ and $B_k = B(\theta_k)$ are known exactly.
When predicting the system behavior at step $k+j$ using the nominal parameter values $\bar{\theta}_{k+j}$, the system may be affected by both an additive disturbance $w_{k+j} \in \mathbb{W}$ and a parametric uncertainty $\Delta\theta_{k+j} \in \Theta_\Delta$, for all $j=1,2,\dots,N$.
%Considering  future steps $k+j$, the system may be affected by both an additive disturbance $w_{k+j} \in \mathbb{W}$ and the parametric uncertainty $\Delta\theta_{k+j} \in \Theta_\Delta$  due to the use of $\bar{\theta}_{k+j}$. 
To handle this, we lump the effects of the additive disturbance and the parametric uncertainty into a single augmented additive disturbance, $d^\theta$. This leads to the following system representation,
\begin{equation}\label{Equ:LPVSystemWAdditiveDis+ParametricUns}
x_{k+j+1} = \bar{A}_{k+j} x_{k+j} + \bar{B}_{k+j} u_{k+j} + d^\theta_{k+j}\ec
\end{equation}
where $d^\theta_{k+j} \in \BD_\theta$ is defined as
\[
d^\theta_{k+j}\triangleq \left(\sum_{i=1}^{p}\Delta\theta_{k}[i]A^{(i)}\right)x_k+\left(\sum_{i=1}^{p}\Delta\theta_{k}[i]A^{(i)}\right)u_k+w_k\ec
\]
with
\begin{equation}
    \BD_\theta\triangleq \conv{\left.\left(\sum_{i=1}^{p}\Delta\theta_{k}[i]A^{(i)}\right)x_k+\left(\sum_{i=1}^{p}\Delta\theta_{k}[i]A^{(i)}\right)u_k +w_k \right| \Delta\theta_k\in\Theta_\Delta\ec w_k\in\BW\ec x_k\in\BsX\ec u_k\in\BU}\ef
\end{equation}
This is again a polytopic set, based on Assumption~\ref{The:Constraints}.

 % is a modified disturbance that accounts for both the original additive disturbance and the uncertainty arising from the estimation error of $\theta_{k+j}$.

The formulation for the Robust Case has been previously considered\cite{Rawlings:2013}. However, it may lead to overly conservative behavior in the LPV Case. Other alternative methods exist for converting parametric uncertainties into additive disturbances using worst-case bounds, which can reduce conservatism depending on the choice of other design parameters in the MPC formulation, such as terminal constraints \cite{Bujarbaruah:2021}. In contrast, considering the above formulation in the LPV Case does not necessarily result in overly conservative behavior, due to Assumption~L.\ref{The:Measurement-theta} \cite{Abbas:2019}. In general, considering this simple approach is practical and facilitates the subsequent formulation of the proposed DRMPC schemes.

%%%------------------------------------------------------------------------------
\section{Deadbeat Robust Model Predictive Control (DRMPC)}

\subsection{Deadbeat Disturbance Feedback Policy}\label{Sec:DeadbeatPolicy}

A slightly different formulation of DRMPC is proposed in this paper, compared to the use of \emph{deadbeat input sequences} for additive disturbances \cite{Schildbach:2025}. Instead, a \emph{deadbeat (disturbance feedback) policy} for pre-stabilization is introduced % for the LPV case, 
as explained below.

Consider the nominal system of \eqref{Equ:DTSystem},
\begin{equation}\label{Equ:NominalSystem}
  x_{k+1}=\bar{A}x_{k}+\bar{B}u_{k}\ec
\end{equation}
starting from an initial condition $x_{0}=d_{0}$ that corresponds to a single disturbance $d_{0}\in\mathbb{D}$, and the following definition.

\begin{definition}[deadbeat horizon]\label{The:DeadbeatHorizon}
The \emph{deadbeat horizon} $M$ \cite{Schildbach:2025} is the smallest number such that the matrix $\mathcal{P}_{M}\in\BR^{n\times nm}$,
\begin{equation}\label{Equ:ControllabilityMatrix}
  \mathcal{P}_{M}\eqdef
  \Bigl[
  \begin{array}{c|c|c|c|c}
     \;\bar{A}^{M-1}\bar{B}\; & \;\bar{A}^{M-2}\bar{B}_{0}\; & \;\cdots\; & \;\bar{A}\bar{B}\; & \;\bar{B}\; 
  \end{array}
  \Bigr]\ec
\end{equation}
has full rank $n$.
\end{definition}
By virtue of Assumption R.\ref{The:Controllability}, $M$ is well defined and satisfies $M\leq n$, because $\mathcal{P}_{M}$ is equivalent to the \emph{controllability matrix} of \eqref{Equ:NominalSystem}.  This leads to the following result.

\begin{lemma}[deadbeat disturbance feedback policy]\label{The:DeadbeatPolicy}
There exists a \emph{deadbeat (disturbance feedback) policy}, consisting of $M$ linear feedback gains $\bar{K}_{0},\bar{K}_{1},\dots,\bar{K}_{M-1}\in\BR^{m\times n}$ with
\begin{equation}\label{Equ:DeadbeatPolicy}
  \bar{u}_{0}=\bar{K}_{0}d_{0}\ec\qquad\bar{u}_{1}=\bar{K}_{1}d_{0}\ec\qquad\dots\ec\qquad\bar{u}_{M-1}=\bar{K}_{M-1}d_{0}\ec
\end{equation}
such that $x_{M}=0$ for any initial condition $d_{0}$.
\end{lemma}

\emph{\textbf{Proof:}}
First, observe that the first $M$ states of the nominal system \eqref{Equ:NominalSystem} can be expressed as
\begin{subequations}\label{Equ:FirstMStates}\begin{align}
  x_{1} &= \bar{A}d_{0}+\bar{B}u_{0}\ec\\
  x_{2} &= \bar{A}^{2}d_{0}+\bar{A}\bar{B}u_{0}+\bar{B}u_{1}\ec\\
  &\cdots\nonumber\\
  x_{M-1} &= \bar{A}^{M-1}d_{0}+\bar{A}^{M-2}\bar{B}u_{0}+\dots+\bar{A}\bar{B}u_{M-3}+\bar{B}u_{M-2}\ec\\
  x_{M} &= \bar{A}^{M}d_{0}+\bar{A}^{M-1}\bar{B}u_{0}+\bar{A}^{M-2}\bar{B}u_{1}+\dots+\bar{A}\bar{B}u_{M-2}+\bar{B}u_{M-1}\;.
\end{align}\end{subequations}
Since the matrix $\mathcal{P}_{M}$ in \eqref{Equ:ControllabilityMatrix} has full rank, there exists a sequence of deadbeat inputs $\bar{u}_{0},\bar{u}_{1},\dots,\bar{u}_{M-1}$ that force $x_{M}=0$. Using $x_{M}\eqcon 0$, (\ref{Equ:FirstMStates}d) becomes
\begin{equation}\label{Equ:DeadbeatCondition1}
  -\bar{A}^{M}d_{0} \eqcon \bar{A}^{M-1}\bar{B}\bar{u}_{0}+\bar{A}^{M-2}\bar{B}u_{1}+\dots+\bar{A}\bar{B}\bar{u}_{M-2}+\bar{B}\bar{u}_{M-1}\;.
\end{equation}
To verify that the deadbeat policy is indeed well-defined, substitute \eqref{Equ:DeadbeatPolicy} into \eqref{Equ:DeadbeatCondition1} to obtain
\begin{equation}\label{Equ:DeadbeatCondition2}
  -\bar{A}^{M}d_{0} = 
  \Bigl[
  \begin{array}{c|c|c|c|c}
  \bar{A}^{M-1}\bar{B}\bar{K}_{0} 
  & \bar{A}^{M-2}\bar{B}\bar{K}_{1}
  &\dots
  &\bar{A}\bar{B}\bar{K}_{M-2}
  &\bar{B}\bar{K}_{M-1}
  \end{array}
  \Bigr]d_{0}\;.
\end{equation}
Since the deadbeat policy must work for all $d_{0}\in\BR^{n}$, the matrices multiplying $d_{0}$ in \eqref{Equ:DeadbeatCondition2} on both sides must be equal, thus
\begin{equation}\label{Equ:DeadbeatCondition3}
  -\bar{A}^{M} = 
  \Bigl[
  \begin{array}{c|c|c|c|c}
  \bar{A}^{M-1}\bar{B}
  & \bar{A}^{M-2}\bar{B}
  &\dots
  &\bar{A}\bar{B}
  &\bar{B}
  \end{array}
  \Bigr]
  \left[
  \begin{array}{c}
  \bar{K}_{0} \\\hline
  \bar{K}_{1}\\\hline
  \dots\\\hline
  \bar{K}_{M-2}\\\hline
  \bar{K}_{M-1}
  \end{array}
  \right]\;.
\end{equation}
This can be expressed as a linear equation system to solve for the linear feeback gains $\bar{K}_{0},\bar{K}_{1},\dots,\bar{K}_{M-1}$:
\begin{equation}\label{Equ:DeadbeatCondition4}
  \left[
  \begin{array}{ccccc}
    \mathcal{P}_{M} & 0 & \cdots & 0 & 0\\
    0 & \mathcal{P}_{M} & \cdots & 0 & 0\\
    \vdots & & \ddots & & \vdots\\
    0 & 0 & \cdots & \mathcal{P}_{M} & 0\\
    0 & 0 & \cdots & 0 & \mathcal{P}_{M}
  \end{array}\right]
  \vct{\left[
  \begin{array}{c}
  \bar{K}_{0} \\\hline
  \bar{K}_{1}\\\hline
  \dots\\\hline
  \bar{K}_{M-2}\\\hline
  \bar{K}_{M-1}
  \end{array}
  \right]}
  =
  -\vct{\bar{A}^{M}}\;.
\end{equation}
The claim follows from the fact that $\mathcal{P}_{M}$ has full rank $n$, according to Definition \ref{The:DeadbeatHorizon}, and hence the matrix on the left-hand side of \eqref{Equ:DeadbeatCondition4} has full rank $Mn$.\hfill$\square$

\begin{remark}
  Alternatively, it is possible to calculate a \emph{time-invariant} disturbance feedback gain $\bar{K}\in\BR^{m}$ by pole placement. Placing all $n$ poles at $0$ by an appropriate choice of $\bar{K}$ renders the matrix $(\bar{A}+\bar{B}\bar{K})$ nilpotent; hence the deadbeat property follows. However, this only works for the choice of $M=n$.
\end{remark}

Note that Lemma \ref{The:DeadbeatPolicy} not only proves the existence of the deadbeat policy. It also provides a constructive approach for computing the corresponding feedback gains $\bar{K}_{0},\bar{K}_{1},\dots,\bar{K}_{M-1}$, namely by solving the linear equation system \eqref{Equ:DeadbeatCondition4}. Besides the existence of a deadbeat policy, however, there is no guarantee for its uniqueness.

\begin{remark}\label{The:ChoiceOfM}
  (a) If the deadbeat policy is not unique, then some additional cost function (e.g., least squares) has to be defined when solving for $\bar{K}_{0},\bar{K}_{1},\dots,\bar{K}_{M-1}$, subject to the equality constraint \eqref{Equ:DeadbeatCondition4}.
  (b) The deadbeat horizon $M$ may be deliberately chosen as a higher number than the smallest number for which $\mathcal{P}_{M}$ has full rank $n$. This opens (additional) degrees of freedom in the choice of the deadbeat policy.
\end{remark}

\subsection{Finite-time Optimal Control Problem}\label{Sec:FTOCP}

For the setup of the Finite-time Optimal Control Problem (FTOCP), it is assumed that a sequence of deadbeat feedback gains $\bar{K}_{0},\bar{K}_{1},\dots,\bar{K}_{M-1}$ have been pre-computed offline, as described in Section \ref{Sec:DeadbeatPolicy}. Moreover, a prediction horizon $N\geq M$ must be selected, typically $N\gg M$.

Similar to \cite{Schildbach:2025}, the setup of the LPV-DRMPC problem is described only for step $k=0$, as the generalization to subsequent steps $k=1,2,\dots$ is straightforward, see Section~\ref{sec:LPV-DRMPC}. Based on Lemma \ref{The:DeadbeatPolicy}, the following affine disturbance feedback policy should be applied in order to account for the deadbeat of each disturbance over the prediction horizon:
\begin{subequations}\label{Equ:DistCompensation}\begin{align}
  u_{0}&=u_{0|0}\ec\\
  u_{1}&=u_{1|0}\,\,+\,\,\bar{K}_{0}d_{0}\ec\\
  u_{2}&=u_{2|0}\,\,+\,\,\bar{K}_{1}d_{0}\,\,+\,\,\bar{K}_{0}d_{1}\ec\\
  u_{3}&=u_{3|0}\,\,+\,\,\bar{K}_{2}d_{0}\,\,+\,\,\bar{K}_{1}d_{1}\,\,+\,\,\bar{K}_{0}d_{2}\ec\\
  &\;\;\vdots\nonumber\\
  u_{M}&=u_{M|0}\,\,+\,\,\bar{K}_{M-1}d_{0}\,\,+\,\,\hdots\,\,+\,\,\bar{K}_{1}d_{M-2}\,\,+\,\,\bar{K}_{0}d_{M-1}\ec\hspace*{0.9cm}\phantom{a}
\end{align}\end{subequations}
\noindent up to the deadbeat horizon, and beyond that
\addtocounter{equation}{-1}
\begin{subequations}\setcounter{equation}{5}\begin{equation}
  u_{j}=u_{j|0}\,\,+\,\,\sum_{i=0}^{M-1}\bar{K}_{i}d_{j-i-1}\qquad\fa M\leq j\leq N-1\ef\hspace*{2.1cm}\phantom{a}
\end{equation}\end{subequations}
The notation $u_{k+j|k}$ is used for the nominal control input in step $k+j$, as planned in step $k$, and $x_{k+j|k}$ for the corresponding predicted nominal state, where $k,j\in\BZ_{0+}$. This scheme suggests the following tightening of the input constraints for the nominal control inputs:
\begin{subequations}\label{Equ:TightenedIC}\begin{align}
  u_{0|0}&\in\BU\ec\\
  u_{1|0}&\in\BU\ominus\bar{K}_{0}\BD\ec\\
  u_{2|0}&\in\BU\ominus\bar{K}_{0}\BD\ominus\bar{K}_{1}\BD\ec\\
  u_{3|0}&\in\BU\ominus\bar{K}_{0}\BD\ominus\bar{K}_{1}\BD\ominus\bar{K}_{2}\BD\ec\\
  &\;\;\vdots\nonumber\\
  u_{M|0}&\in\BU\ominus\bar{K}_{0}\BD\ominus\bar{K}_{1}\BD\ominus\hdots\ominus\bar{K}_{M-1}\BD\ec
\end{align}\end{subequations}
\noindent up to the deadbeat horizon, and beyond that
\addtocounter{equation}{-1}
\begin{subequations}\setcounter{equation}{5}\begin{equation}
  \hspace*{3.1cm}u_{j|0}\in\BU\ominus\bar{K}_{0}\BD\ominus\bar{K}_{1}\BD\ominus\cdots\ominus\bar{K}_{M-1}\BD\qquad\fa M\leq j\leq N-1\ef
\end{equation}\end{subequations}
For the tightening of the state constraints, define the auxiliary matrices
\begin{subequations}\label{Equ:AuxMatrices}\begin{align}
  \bar{\Phi}_{0}&\eqdef\bar{A}+\bar{B}\bar{K}_{0}\ec\\
  \bar{\Phi}_{1}&\eqdef\bar{A}^2+\bar{A}\bar{B}\bar{K}_{0}+\bar{B}\bar{K}_{1}\ec\\
  \bar{\Phi}_{2}&\eqdef\bar{A}^3+\bar{A}^2\bar{B}\bar{K}_{0}+\bar{A}\bar{B}\bar{K}_{1}+\bar{B}\bar{K}_{2}\ec\\
  &\;\;\vdots\nonumber\\
\bar{\Phi}_{M-2}&\eqdef\bar{A}^{M-1}+\bar{A}^{M-2}\bar{B}\bar{K}_{0}+\bar{A}^{M-3}\bar{B}\bar{K}_{1}+\dots+\bar{A}\bar{B}\bar{K}_{M-3}+\bar{B}\bar{K}_{M-2}\ef
\end{align}\end{subequations}
Note that, by construction, the next auxiliary matrix $\bar{\Phi}_{M-1}=0$, due to \eqref{Equ:DeadbeatCondition3}. With the help of these auxiliary matrices, the constraints for the nominal states should be tightened as follows:
\begin{subequations}\label{Equ:TightenedSC}\begin{align}
  x_{1|0}&\in\BsX\ominus\BD\ec\\
  x_{2|0}&\in\BsX\ominus\BD\ominus\bar{\Phi}_{0}\BD\ec\\
  x_{3|0}&\in\BsX\ominus\BD\ominus\bar{\Phi}_{0}\BD\ominus\bar{\Phi}_{1}\BD\ec\\
  &\;\;\vdots\nonumber\\
  x_{M|0}&\in\BsX\ominus\BD\ominus\bar{\Phi}_{0}\BD\ominus\bar{\Phi}_{1}\BD\ominus\hdots\ominus\bar{\Phi}_{M-2}\BD\ec
\end{align}\end{subequations}
\noindent up to the deadbeat horizon, and beyond that
\addtocounter{equation}{-1}
\begin{subequations}\setcounter{equation}{5}\begin{equation}
  \hspace*{2.8cm}x_{j|0}\in\BsX\ominus\BD\ominus\bar{\Phi}_{0}\BD\ominus\bar{\Phi}_{1}\BD\ominus\hdots\ominus\bar{\Phi}_{M-2}\BD\qquad\fa M\leq j\leq N\ef
\end{equation}\end{subequations}

Combining these constraints with a stage cost function, based on a \emph{stage cost} $\ell:\BsX\times\BU\to\BR$, yields the following \emph{DRMPC (optimization) problem}:
\begin{subequations}\label{Equ:DRMPCProblem}\begin{align}
  \min\enspace&\sum_{j=0}^{N-1}\ell\left(u_{j|0},x_{j|0}\right)\\
  \text{s.t.}\enspace &\quad x_{j+1|0}=\bar{A}x_{j|0}+\bar{B}u_{j|0}\qquad\fa\,j\in\{0,\hdots,N-1\}\ec\\
   &\quad x_{0|0}=x_{0}\ec\\
   &\quad\eqref{Equ:TightenedIC}\ec\quad\eqref{Equ:TightenedSC}\ec\\
   &\quad x_{N|0}=0\ef
\end{align}\end{subequations}

%%------------------------------------------------------------------------
\section{DRMPC for LPV Systems (LPV-DRMPC)}\label{sec:LPV-DRMPC}

\subsection{Deadbeat Disturbance Feedback Policy for LPV Systems }\label{Sec:DeadbeatPolicy-LPV}

Consider the nominal LPV system \eqref{Equ:NominalSystemLPV}.
%\begin{equation}\label{Equ:NominalSystemLPV}
%  x_{k+1}=\bar{A}_kx_{k}+\bar{B}_ku_{k}.
%\end{equation}
As in the Robust Case, starting from an initial condition $x_{0}=d_0$, which corresponds to a single additive  disturbance $d_{0}$, consider the following definition as
 the LPV counterpart to the Robust Case.
%due to the error in the given (estaimted) value of $\theta_k$ ({\color{red} should start with $d_1^\theta$ not $d_0^\theta$ as $\theta$ is known at step $k$, I will update later the following drivation}).
\begin{definition}[deadbeat horizon in the LPV Case]\label{The:DeadbeatHorizonLPV}
The \emph{deadbeat horizon} $M$ \cite{Schildbach:2025} is the smallest number such that the matrix $\mathcal{P}_{M}^\theta\in\BR^{n\times nm}$,
\begin{equation}\label{Equ:ControllabilityMatrixLPV-2}
  \mathcal{P}_{M}^\theta\eqdef
  \Bigl[
  \begin{array}{c|c|c|c|c}
     \;\prod_{i=1}^{M-1}\bar{A}_i\bar{B}_0\; & \;\prod_{i=2}^{M-1}\bar{A}_i\bar{B}_1\; & \;\cdots\; & \;\bar{A}_{M-1}\bar{B}_{M-2}\; & \;\bar{B}_{M-1}\; 
  \end{array}
  \Bigr]\ec
\end{equation}
has full rank $n$.
\end{definition}
\noindent
According to Assumption~\ref{The:ControllabilityLPV}, $M$ is well defined and satisfies $M\geq n$   due to the structural equivalence between    $\mathcal{P}_{M}^\theta$ and the controllability matrix $\mathcal{C}_{k,r}$ in \eqref{Equ:ControllabilityMatrixLPV-1}. Therefore, Lemma~\ref{The:DeadbeatPolicy} applies in the LPV Case, as shown in the following proof.

\emph{\textbf{Proof:}} The first $M$ states of the nominal system \eqref{Equ:NominalSystemLPV} can be expressed as
\begin{subequations}\label{Equ:FirstMStatesLPV}
\begin{align}
  x_{1} &= \bar{A}_0d_{0}+\bar{B}_0u_{0}\ec\\
  x_{2} &= \bar{A}_1\bar{A}_0d_{0}+\bar{A}_1\bar{B}_0u_{0}+\bar{B}_1u_{1}\ec\\
  &\vdots\nonumber\\
  x_{M-1} &= \prod_{i=0}^{M-2}\bar{A}_id_{0}+\prod_{i=1}^{M-2}\bar{A}_i\bar{B}_0u_{0}+\dots+\bar{A}_{M-2}\bar{B}_{M-3}u_{M-3}+\bar{B}_{M-2}u_{M-2}\ec\\
  x_{M} &= \prod_{i=0}^{M-1}\bar{A}_i d_{0}+\prod_{i=1}^{M-1}\bar{A}_i\bar{B}_0u_{0}+\prod_{i=2}^{M-1}\bar{A}_i\bar{B}_1u_{1}+\dots+\bar{A}_{M-1}\bar{B}_{M-2}u_{M-2}+\bar{B}_{M-1}u_{M-1}\;.
\end{align}
\end{subequations}
Since the matrix $\mathcal{P}_{M}^\theta$ in \eqref{Equ:ControllabilityMatrixLPV-2} has full rank, there exists a sequence of deadbeat inputs $\bar{u}_{0},\bar{u}_{1},\dots,\bar{u}_{M-1}$ that force $x_{M}=0$. Using $x_{M}\eqcon 0$, (\ref{Equ:FirstMStatesLPV}d) becomes
\begin{equation}\label{Equ:DeadbeatConditionLPV1}
  -\prod_{i=0}^{M-1}\bar{A}_i d_0\eqcon \prod_{i=1}^{M-1}\bar{A}_i\bar{B}_0u_{0}+\prod_{i=2}^{M-1}\bar{A}_i\bar{B}_1u_{1}+\dots+\bar{A}_{M-1}\bar{B}_{M-2}u_{M-2}+\bar{B}_{M-1}u_{M-1}\;.
\end{equation}
Since the deadbeat policy must work for all $d_{0}\in\BR^{n}$, the matrices multiplying $d_{0}$ in \eqref{Equ:DeadbeatCondition2} on both sides must be equal, thus
\begin{equation}\label{Equ:DeadbeatCondition3LPV}
  -\prod_{i=0}^{M-1}\bar{A}_i = 
  \Bigl[
  \begin{array}{c|c|c|c|c}
     \;\prod_{i=1}^{M-1}\bar{A}_i\bar{B}_0\; & \;\prod_{i=2}^{M-1}\bar{A}_i\bar{B}_1\; & \;\cdots\; & \;\bar{A}_{M-1}\bar{B}_{M-2}\; & \;\bar{B}_{M-1}\; 
  \end{array}  \Bigr]
  \left[
  \begin{array}{c}
  \bar{K}^\theta_{0} \\\hline
  \bar{K}^\theta_{1}\\\hline
  \dots\\\hline
  \bar{K}^\theta_{M-2}\\\hline
  \bar{K}^\theta_{M-1}
  \end{array}
  \right]\;.
\end{equation}
This can be expressed as a linear system of equations, similarly to \eqref{Equ:DeadbeatCondition4}, to solve for the linear feedback gains $\bar{K}^\theta_{0}, \bar{K}^\theta_{1}, \dots, \bar{K}^\theta_{M-1}$, using the fact that $\mathcal{P}_{M}^\theta$ has full rank $n$ according to Definition~\ref{The:DeadbeatHorizonLPV}.\hfill$\square$

\subsection{Finite-time Optimal Control Problem for  LPV Systems}\label{Sec:FTOCP}
Based on Lemma \ref{The:DeadbeatPolicy} for the LPV Case, the affine disturbance feedback policy is applied   in order to account for the deadbeat of each disturbance over the prediction horizon. Considering Assumption~L.\ref{The:Measurement-theta} and the conversion to additive disturbance properties in Section~\ref{Sec:AdditiveUncertainty}, at step $i = 0$ of the prediction horizon, only the additive disturbance $w_0 \in \mathbb{W}$ is considered, and the system is represented by \eqref{Equ:LPVSystemWAdditiveDisOnly}. At steps $i \geq 1$, the augmented additive disturbance $d_i^\theta \in \mathbb{D}_\theta$ is considered along with the nominal part of \eqref{Equ:LPVSystemWAdditiveDis+ParametricUns}. Therefore,  the following affine disturbance feedback policy should be applied using the gains obtained from \eqref{Equ:DeadbeatCondition3}:
\begin{subequations}\label{Equ:DistCompensation-LPV}\begin{align}
  u_{0}&=u_{0|0}\ec\\
  u_{1}&=u_{1|0}\,\,+\,\,\bar{K}^\theta_{0}w_{0}\ec\\
  u_{2}&=u_{2|0}\,\,+\,\,\bar{K}^\theta_{1}w_{0}\,\,+\,\,\bar{K}^\theta_{0}d^\theta_{1}\ec\\
  u_{3}&=u_{3|0}\,\,+\,\,\bar{K}^\theta_{2}w_{0}\,\,+\,\,\bar{K}^\theta_{1}d^\theta_{1}\,\,+\,\,\bar{K}^\theta_{0}d^\theta_{2}\ec\\
  &\;\;\vdots\nonumber\\
  u_{M}&=u_{M|0}\,\,+\,\,\bar{K}^\theta_{M-1}w_{0}\,\,+\,\,\bar{K}^\theta_{M-2}d^\theta_{1}\,\,+\,\,\hdots\,\,+\,\,\bar{K}^\theta_{1}d^\theta_{M-2}\,\,+\,\,\bar{K}^\theta_{0}d^\theta_{M-1}\ec\hspace*{0.9cm}\phantom{a}
\end{align}\end{subequations}
\noindent up to the deadbeat horizon, and beyond that
\addtocounter{equation}{-1}
\begin{subequations}\setcounter{equation}{5}\begin{equation}
  u_{j}=u_{j|0}\,\,+\,\,\,\bar{K}^\theta_{M-1}w_{0}\,+\,\,\sum_{i=0}^{M-2}\bar{K}^\theta_{i}d^\theta_{j-i-1}\qquad\fa M\leq j\leq N-1\ef\hspace*{2.1cm}\phantom{a}
\end{equation}\end{subequations}
Therefore, the following tightening of the input constraints for the nominal control inputs is applied:
\begin{subequations}\label{Equ:TightenedIC-LPV}\begin{align}
  u_{0|0}&\in\BU\ec\\
  u_{1|0}&\in\BU\ominus\bar{K}^\theta_{0}\BW\ec\\  u_{2|0}&\in\BU\ominus\bar{K}^\theta_{1}\BW\ominus\bar{K}^\theta_{0}\BDt\ec\\ u_{3|0}&\in\BU\ominus\bar{K}^\theta_{2}\BW\ominus\bar{K}^\theta_{1}\BDt\ominus\bar{K}^\theta_{0}\BDt\ec\\
  &\;\;\vdots\nonumber\\  u_{M|0}&\in\BU\ominus\bar{K}^\theta_{M-1}\BW\ominus\bar{K}^\theta_{M-2}\BDt\ominus\hdots\ominus\bar{K}^\theta_{0}\BDt\ec
\end{align}\end{subequations}
\noindent up to the deadbeat horizon, and beyond that
\addtocounter{equation}{-1}
\begin{subequations}\setcounter{equation}{5}\begin{equation}
  \hspace*{3.1cm}u_{j|0}\in\BU\ominus\bar{K}^\theta_{M-1}\BW\ominus\bar{K}^\theta_{M-2}\BDt\ominus\cdots\ominus\bar{K}^\theta_{0}\BDt\qquad\fa M\leq j\leq N-1\ef
\end{equation}\end{subequations}

The corresponding constraints for the nominal  states in the LPV Case should be tightened as follows:
\begin{subequations}\label{Equ:TightenedSCLPV}\begin{align}
  x_{1|0}&\in\BsX\ominus\BW\ec\\
  x_{2|0}&\in\BsX\ominus\bar{\Phi}^\theta_0\BW\ominus\BDt\ec\\
x_{3|0}&\in\BsX\ominus\bar{\Phi}^\theta_1\BW\ominus\bar{\Phi}^\theta_0\BDt\ominus\BDt\ec\\
  &\;\;\vdots\nonumber\\  x_{M|0}&\in\BsX\ominus\bar{\Phi}^\theta_{M-2}\BW\ominus\bar{\Phi}^\theta_{M-3}\BDt\ominus\hdots\ominus\bar{\Phi}^\theta_1\BDt\ominus\bar{\Phi}^\theta_0\BDt\ominus\BDt\ec
\end{align}\end{subequations}
\noindent up to the deadbeat horizon, and beyond that
\addtocounter{equation}{-1}
\begin{subequations}\setcounter{equation}{5}\begin{equation}
  \hspace*{2.8cm}x_{j|0}\in\BsX\ominus\bar{\Phi}^\theta_{M-2}\BW\ominus\bar{\Phi}^\theta_{M-3}\BDt\ominus\hdots\ominus\bar{\Phi}^\theta_1\BDt\ominus\bar{\Phi}^\theta_0\BDt\ominus\BDt\qquad\fa M\leq j\leq N\ec
\end{equation}\end{subequations}
where the auxiliary matrices are given as follows:
\begin{subequations}\label{Equ:AuxMatrices-LPV}\begin{align}
  \bar{\Phi}^\theta_{0}&\eqdef\bar{A}_0+\bar{B}_0\bar{K}_{0}\ec\\  \bar{\Phi}^\theta_{1}&\eqdef\bar{A}_1\bar{A}_0+\bar{A}_1\bar{B}_0\bar{K}_{0}+\bar{B}_1\bar{K}_{1}\ec\\  \bar{\Phi}^\theta_{2}&\eqdef\bar{A}_2\bar{A}_1\bar{A}_0+\bar{A}_2\bar{A}_1\bar{B}_0\bar{K}_{0}+\bar{A}_2\bar{B}_1\bar{K}_{1}+\bar{B}_2\bar{K}_{2}\ec\\
  &\;\;\vdots\nonumber\\
\bar{\Phi}^\theta_{M-2}&\eqdef\prod_{i=0}^{M-1}\bar{A}_i+\prod_{i=1}^{M-1}\bar{A}_i\bar{B}_{0}\bar{K}_{0}+\prod_{i=2}^{M-1}\bar{A}_i\bar{B}_{1}\bar{K}_{1}+\dots+\bar{A}_{M-1}\bar{B}_{M-3}\bar{K}_{M-3}+\bar{B}_{M-2}\bar{K}_{M-2}\ef
\end{align}\end{subequations}
%with  $\bar{\Phi}^\theta_{M-1}=0$ according to \eqref{Equ:DeadbeatCondition3LPV}.
Combining these constraints with a stage cost function, based on a \emph{stage cost} $\ell:\BsX\times\BU\to\BR$, yields the following \emph{LPV-DRMPC  optimization problem}:
\begin{subequations}\label{Equ:DRMPCP-LPVroblem}\begin{align}
  \min\enspace&\sum_{j=0}^{N-1}\ell\left(u_{j|0},x_{j|0}\right)\\
  \text{s.t.}\enspace &\quad x_{j+1|0}=\bar{A}_jx_{j|0}+\bar{B}_ju_{j|0}\qquad\fa\,j\in\{0,\hdots,N-1\}\ec\\
   &\quad x_{0|0}=x_{0}\ec\\
   &\quad\eqref{Equ:TightenedIC-LPV}\ec\quad\eqref{Equ:TightenedSCLPV}\ec\\
   &\quad x_{N|0}=0\ef
\end{align}\end{subequations}

%-----------------------------------------------------------------------------------------------------------------------
\section{System-Theoretical Analysis}\label{Sec:SystemTheory}
%-----------------------------------------------------------------------------------------------------------------------

\subsection{The Robust Case}\label{Sec:SystemTheory-Robustcase}

%\section{DRMPC for Linear Systems (DRMPC)}

As a preliminary step, consider the special case of a linear system with additive disturbances \eqref{Equ:DTSystem}; i.e., the parameter $\theta_{k}$ is time-invariant and known, and the parameter set $\Theta$ consists only of a single point. As shown in this section, in this case, the DRMPC comes with all relevant system-theoretic guarantees, including recursive feasibility and stability. The line of argument is closely related to the previous work \cite{Schildbach:2025}.
%-----------------------------------------------------------------------------------------------------------------------
\subsubsection{Recursive Feasibility}\label{Sec:LinearFeasibility}

Under the DRMPC regime, in each time step $k=0,1,\hdots$ the first element of the optimal input sequence from \eqref{Equ:DRMPCProblem}, $u_{k}=u_{k|k}^{\star}$, is applied to the system.

\begin{theorem}[recursive feasibility for the Robust Case]\label{The:LinearRecFeasibility}
  Assume system \eqref{Equ:DTSystem} is linear with additive disturbances, i.e., $d=d_0$ for some $d_0\in\BD$.  %$\Theta=\bar{\theta}$ for some $\bar{\theta}\in\BR^{p}$.
  If the DRMPC problem \eqref{Equ:DRMPCProblem} is feasible at $k=0$ for $x=0$, % $\bar{x}_{0}$, 
  it remains feasible for all future states $x_{1},x_{2},\hdots$ of system \eqref{Equ:DTSystem} under the DRMPC regime.
\end{theorem}

\emph{\textbf{Proof:}} Consider the solution to \eqref{Equ:DRMPCProblem} in step $k=0$. Let $u_{0|0}^{\star},\hdots,u_{N-1|0}^{\star}$ be the optimal inputs and $x_{0|0}^{\star},\hdots,x_{N|0}^{\star}$ be the optimal states. They satisfy the input constraints \eqref{Equ:TightenedIC} and state constraints \eqref{Equ:TightenedSC}, respectively. By the terminal constraint (\ref{Equ:DRMPCProblem}e), $x_{N|0}^{\star}=0$. 

For recursive feasibility, it must be shown that in time step $k=1$ there exists a candidate input sequence $u_{1|1},\hdots,u_{N|1}$ and a corresponding state sequence $x_{1|1},\hdots,x_{N+1|1}$ that satisfy the input, state, and terminal constraints, for any admissible disturbance $d_{0}$. By the presence of a disturbance,
\begin{equation*}
  x_{1|1}=x_{1|0}^{\star}+d_{0}\ec\quad\text{for some}\;d_{0}\in\BD\ef
\end{equation*}

Consider the following candidate input sequence for step $k=1$,
\begin{align*}
  u_{1|1}&=u_{1|0}^{\star}+\bar{K}_{0}d_{0}\ec\\
  u_{2|1}&=u_{2|0}^{\star}+\bar{K}_{1}d_{0}\ec\\
  &\;\;\vdots\\
  u_{M|1}&=u_{M|0}^{\star}+\bar{K}_{M-1}d_{0}\ec\\
  u_{M+1|1}&=u_{M+1|0}^{\star}\ec\\
  &\;\;\vdots\\
  u_{N-1|1}&=u_{N-1|0}^{\star}\ec\\
  u_{N|1}&=0\ef
\end{align*}
This candidate input sequence $u_{1|1},\dots,u_{N|1}$ satisfies the input constraints (\ref{Equ:TightenedIC}a-f), given that the input sequence $u_{0|0}^{\star},\hdots,u_{N-1|0}^{\star}$ satisfies (\ref{Equ:TightenedIC}a-f). The corresponding state sequence is given by
\begin{align*}
  x_{1|1}&=x_{1|0}^{\star}+d_{0}\ec\\
  x_{2|1}&=x_{2|0}^{\star}+\bar{\Phi}_{0}d_{0}\ec\\
  x_{3|1}&=x_{3|0}^{\star}+\bar{\Phi}_{1}d_{0}\ec\\
  &\;\;\vdots\\
  x_{M|1}&=x_{M|0}^{\star}+\bar{\Phi}_{M-2}d_{0}\ec\\
  x_{M+1|1}&=x_{M+1|0}^{\star}\ec\\
  x_{M+2|1}&=x_{M+2|0}^{\star}\ec\\
  &\;\;\vdots\\
  x_{N|1}&=x_{N|0}^{\star}=0\ec\\
  x_{N+1|1}&=0\ef
\end{align*}
\noindent It satisfies all the state constraints (\ref{Equ:TightenedSC}a-f) and the terminal constraint (\ref{Equ:DRMPCProblem}e). The last step follows from $x_{N|0}^{\star}=0$ and the choice of $u_{N|1}=0$.\hfill$\square$
%-----------------------------------------------------------------------------------------------------------------------
\subsubsection{Input-to-State Stability}\label{Sec:Stability}

Let $\CX_{N}\subseteq\BsX$ be the \emph{feasible set}, i.e., the set of initial conditions $x_{0}$ for which problem \eqref{Equ:DRMPCProblem} is feasible. The \emph{state feedback control law} $\kappa_{N}:\CX_{N}\to\BU$ of the DRMPC regime is defined as $\kappa_{N}(x_{k})=u_{k|k}^{\star}$. From Theorem \ref{The:LinearRecFeasibility}, this map is well-defined.

\begin{definition}[robust positive invariant set]
  Consider system \eqref{Equ:DTSystem} with any state feedback control law $\kappa:\BsX\to\BU$. A set $\CP\subseteq\BsX$ is called a \textbf{robust positive invariant (RPI) set} if for all $x\in\CP$ it holds that
  \begin{equation*}
    \kappa(x)\in\BU\quad\text{and}\quad Ax+B\kappa(x)\oplus\BD\subseteq\CP\;.
  \end{equation*}
\end{definition}

\begin{corollary}[invariance of the feasible set]\label{The:InvFeasSet}
  The feasible set $\CX_{N}$ is a RPI set for system \eqref{Equ:DTSystem} under the DRMPC regime.
\end{corollary}

\emph{\textbf{Proof:}} This is immediate from Theorem \ref{The:LinearRecFeasibility}.\hfill$\square$\vspace*{0.2cm}

Computing the feasible set of the DRMPC problem \eqref{Equ:DRMPCProblem} hence leads to an RPI set. On top of this invariance property, it can be shown that the DRMPC stabilizes system \eqref{Equ:DTSystem}. Due to the presence of disturbances, the appropriate stability concept to be used is \emph{input-to-state stability (ISS)} \cite{JiangWang:2001,LimonEtAl:2009,RaMaDi:2018}. It requires the following standard assumption on the stage cost function.

\begin{assumption}[stage cost]\label{Ass:StageCost}
  For the stage cost function $\ell:\BsX\times\BU\rightarrow\BR_{0+}$ there exist two $\textrm{K}_{\infty}$-functions $\varepsilon_{1}$ and $\varepsilon_{2}$ with
  \begin{equation}\label{Equ:LyaFunction1}
    \varepsilon_{1}(\|x\|)\leq \ell\left(x,u\right)\leq \varepsilon_{2}(\|x\|)\quad\fa x\in\BsX\;,\;u\in\BU\;.
  \end{equation}
\end{assumption}

\begin{definition}[input-to-state stability]\label{Def:Stability}
    System \eqref{Equ:DTSystem} under some state feedback law $\kappa:\CP\to\BU$ is \textbf{input-to-state stable (ISS)} on some RPI set $\CP\subseteq\BsX$ if there exist a $\text{KL}$-function $\beta$ and a $\text{K}$-function $\gamma$ such that
    \begin{equation*}
        \|x_{k}\|\leq\beta\bigl(\|x_{0}\|,k\bigr)+\gamma\bigl(\max_{j<k}\|d_{j}\|\bigr) \quad\fa k\in\BZ_{0+}\;,
    \end{equation*}
    for any $x_{0}\in\CP$ and any disturbances $d_{0},d_{1},\hdots\in\BD$.
\end{definition}

Definition \ref{Def:Stability} reduces to \emph{asymptotic stability} if all disturbances vanish beyond some step $k$, i.e., $d_{j}=0$ for all $j\geq k$. If the disturbances are non-zero, the states $x_{k}$ are bounded by norm balls, whose size increases with the maximum value of $\|d_{j}\|$ observed in the past $j<k$.

\begin{definition}[ISS Lyapunov function]\label{Def:LyaFunction}
  $V:\CP\rightarrow\BR_{0+}$ is called an \textbf{ISS Lyapunov function} for system \eqref{Equ:DTSystem} under the state feedback law $\kappa:\CP\to\BU$, if the following conditions are satisfied:\\
    (a) There exist two $\textrm{K}_{\infty}$-functions $\alpha_{1}$ and $\alpha_{2}$ such that
    \begin{equation}\label{Equ:LyaFunction1}
        \alpha_{1}(\|x\|)\leq V\left(x\right)\leq \alpha_{2}(\|x\|)\quad\fa x\in\CP\;.
    \end{equation}
    (b) There exists a $\textrm{K}_{\infty}$-function $\alpha_{3}$ and $\text{K}$-function $\sigma$ with
    \begin{equation}\label{Equ:LyaFunction2}
        V\left(x_{k+1}\right)\leq V\left(x_{k}\right)-\alpha_{3}(\|x_{k}\|)+\sigma(\|d_{k}\|)
    \end{equation}
    for all $k\in\BZ_{0+}$, any $x_{0}\in\CP$, and any $d_{k}\in\BD$.
\end{definition}

Definition \ref{Def:LyaFunction} reduces to a \emph{Lyapunov function} if the term $\sigma(\|d_{k}\|)$ is removed in \eqref{Equ:LyaFunction2}.

\begin{lemma}[Lyapunov stability]\label{The:LyapunovStability}
  If there exists a ISS Lyapunov function $V:\CP\rightarrow\BR_{0+}$ for system \eqref{Equ:DTSystem} under some feedback law $\kappa:\BsX\to\BU$ and on some RPI set $\CP\subseteq\BsX$, then the closed-loop system is ISS.
\end{lemma}

\emph{\textbf{Proof:}} See \cite{JiangWang:2001}.

\begin{theorem}[input-to-state stability]\label{The:Stability}
  The %closed-loop
  system \eqref{Equ:DTSystem} under the DRMPC regime $\kappa_{N}:\CX_{N}\to\BU$ in closed loop is ISS on $\CX_{N}$.
\end{theorem}

\emph{\textbf{Proof:}} By Lemma \ref{The:LyapunovStability}, it suffices to show that the optimal value function of the DRMPC problem $V^{\star}_{N}:\CX_{N}\to\BR_{0+}$ is an ISS Lyapunov function in the sense of Definition \ref{Def:LyaFunction}.

Definition \ref{Def:LyaFunction}(a) is satisfied by virtue of Assumption \ref{Ass:StageCost}. Definition \ref{Def:LyaFunction}(b) is verified as follows. First, for the case of $d_{k}=0$, the shifted optimal input sequence from the previous step
\begin{align*}
  u_{k+1|k+1}&=u^{\star}_{k+1|k}\;,\\
  &\dots\;,\\
  u_{k+N|k+1}&=u^{\star}_{k+N|k}\;,\\
  u_{k+N+1|k+1}&=0
\end{align*}
is feasible, and so is the corresponding state sequence. By the standard argument in MPC stability proofs, the optimal value function hence decreases by more than the first stage cost,
\begin{equation*}
  V^{\star}_{N}\left(x_{k+1}\right)\leq V^{\star}_{N}\left(x_{k}\right)-\ell\left(x_{k},u_{k}\right)\;.
\end{equation*}
Together with Assumption \ref{Ass:StageCost} this verifies the term $-\alpha_{3}(\|x_{k}\|)$ in \eqref{Equ:LyaFunction2}, for the case of $d_{k}=0$.

Second, for the case of $d_{k}\neq 0$, observe that the optimal value function $V^{\star}_{N}(x_{k+1|k}^{\star}+d_{k})$ is continuous in $d_{k}$, for all $d_{k}\in\mathbb{D}$. Together with the compactness of $\mathbb{D}$, this verifies the term $+\sigma(\|d_{k}\|)$ in \eqref{Equ:LyaFunction2}, for the case of $d_{k}\neq 0$.\hfill$\square$

%-----------------------------------------------------------------------------------------------------------------------
\subsection{The LPV Case}

The corresponding results for the LPV Case follow analogously from those in Section~\ref{Sec:SystemTheory-Robustcase}. For brevity, we only state the key formulations.

%\subsubsection{Recursive Feasibility}\label{Sec:RecursiveFeasibility}

\begin{theorem}[Recursive feasibility for the LPV Case]\label{The:LPVRecFeasibility}
  Assume that system \eqref{Equ:DTSystem} is LPV. Suppose   Assumptions~\ref{The:Measurement-x},~L.\ref{The:Measurement-theta},~L.\ref{The:ControllabilityLPV} and \ref{The:Constraints}  as well as  the considerations in Lemma~\ref{The:DeadbeatPolicy} hold.  If the LPV-DRMPC problem \eqref{Equ:DRMPCP-LPVroblem} is feasible at $k=0$ for ${x}_{0}$, it remains feasible for all future states $x_{1},x_{2},\hdots$ of system \eqref{Equ:DTSystem} under the LPV-DRMPC regime.
\end{theorem}

\emph{\textbf{Proof:}} Following the proof of Theorem~\ref{The:LinearRecFeasibility}, given the optimal inputs $u_{0|0}^{\star},\ldots,u_{N-1|0}^{\star}$ and states $x_{0|0}^{\star},\ldots,x_{N|0}^{\star}$ solving \eqref{Equ:DRMPCP-LPVroblem} at time $k=0$ and satisfying all constraints, we show that at $k=1$, there exists a candidate input sequence $u_{1|1},\ldots,u_{N|1}$ and a corresponding state sequence $x_{1|1},\ldots,x_{N+1|1}$ that satisfy the input, state, and terminal constraints in \eqref{Equ:DRMPCP-LPVroblem} for any admissible disturbance $w_0\in\BW$  yielding $x_{1|1}=x_{1|0}^{\star}+w_0$ and an augmented admissible disturbance in the subsequent step $d_1^\theta\in\BDt$. 
The candidate input sequence $u_{1|1}=u_{1|0}^{\star}+\bar{K}^\theta_{0}w_0$, $u_{2|1}=u_{2|0}^{\star}+\bar{K}^\theta_{1}w_0+\bar{K}^\theta_{0}d_1^\theta$, $\dots$, $u_{M|1}=u_{M|0}^{\star}+\bar{K}^\theta_{M-1}w_0+\bar{K}^\theta_{M-2}d_1^\theta$, $u_{M+1|1}=u_{M+1|0}^{\star}$, $\dots$, $u_{N-1|1}=u_{N-1|0}^{\star}$, $u_{N|1}=0$ satisfies the input constraints \eqref{Equ:TightenedIC-LPV} and leads to the state sequence 
$x_{2|1}=x_{2|0}^{\star}+\bar{\Phi}_0^\theta w_0$, $x_{3|1}=x_{3|0}^{\star}+\bar{\Phi}_1^\theta w_0+\bar{\Phi}_0^\theta d_1^\theta$,  $\dots$, 
$x_{M|1}=x_{M|0}^{\star}+\bar{\Phi}_{M-2}^\theta w_0+\bar{\Phi}_{M-3}^\theta d_1^\theta$, $x_{M+1|1}=x_{M+1|0}^{\star}$,  $\dots$,  $x_{N|1}=x_{N|0}^{\star}$, $x_{N+1|1}=0$, which satisifes all the state constraints \eqref{Equ:TightenedSCLPV}.\hfill$\square$

\begin{theorem}[input-to-state stability for the LPV Case]\label{The:Stability-LPV}
  The LPV  system \eqref{Equ:DTSystem} under the LPV-DRMPC regime $\kappa_{N}:\CX_{N}\to\BU$ in closed loop is ISS on $\CX_{N}$, where $\CX_{N}$ represents the feasible set associated with problem \eqref{Equ:DRMPCP-LPVroblem}.
\end{theorem}
\emph{\textbf{Proof:}} This follows from Lemma \ref{The:LyapunovStability}, Assumption \ref{Ass:StageCost}, and the feasibility results of Theorem~\ref{The:LPVRecFeasibility}, and then, analogously to the proof of Theorem~\ref{The:Stability}, it can be shown that the optimal value function of the LPV-DRMPC problem,  $V^{\star}_{N}:\CX_{N}\to\BR_{0+}$, is an ISS Lyapunov function.\hfill$\square$

%\subsection{Input-to-State Stability}\label{Sec:ISS}

%%-------------------------------------------------------------------
\section{Numerical Example}\label{Sec:NumericalExample}

In this section, the performance of the proposed DRMPC and LPV-DRMPC are demonstrated in a small numerical example. We consider a modification of the system of Fleming et al.\ \cite{FlemingEtAl:2015}, which is rearranged to have a box parameter set with four vertices instead of three, as described in other previous contributions \cite{Abbas:2019,Heydari:2021}.
The affine system matrices, see~\eqref{Equ:Parameter}, are given as follows:
\[    \bar{A}=\begin{bmatrix}
         0.2485 &  -1.0355\\     0.8910  &  0.4065
    \end{bmatrix}\ec \quad
    A^{(1)}=\begin{bmatrix}
          -0.0063  & -0.0938\\          0  &  0.0188
    \end{bmatrix}\ec    \quad
    A^{(2)}=\begin{bmatrix}
          -0.0063  & -0.0938\\          0  &  0.0188
    \end{bmatrix}\ec\quad \bar{B}=\begin{bmatrix}
         0.3190 \\   -1.3080
    \end{bmatrix}\ec   
\]
 and $B^{(1)}=B^{(2)}=0$. The parameter set is defined as $\Theta\triangleq [-1,1]\times [-1,1]\subset\mathbb{R}^2$. The constraint sets are $\mathbb{X}\triangleq [-60,60]\times [-41.7, 41.7]\subset\mathbb{R}^2$ and $\mathbb{U}\triangleq [-12.5, 12.5]\subset\mathbb{R}$. In this example, no additive disturbances are considered, i.e., $\BW = \{0\}$. 
The conversion from parametric to additive uncertainty, as described in Section~\ref{Sec:AdditiveUncertainty}, results in the following set of additive uncertainties, defined as
\[
\BD\triangleq \left\{d\in\mathbb{R}^2   \mid \begin{bmatrix}
    0 &  -0.5387 \\
    0.1940  &  0.9701\\
    0 &   0.5387\\
   -0.1940  & -0.9701
\end{bmatrix} d\leq
\begin{bmatrix}
    0.8425\\     0.1455\\     0.8425\\     0.1455
\end{bmatrix}
\right\}\ec
\]
which is a polytope with four vertices. The stage cost function is assumed to be quadratic,
\begin{equation}
        \ell(u,x)=x^\top Q x+ u^\top R u\ec \qquad\text{with}\;
        Q=\begin{bmatrix}
            1&0\\0&1
        \end{bmatrix}\ec\quad R=1\ef
\end{equation}
The deadbeat horizon is chosen as $M = 3$ and the terminal set is the origin, unless otherwise stated.

The performance of the proposed approaches is evaluated based on the size of their corresponding \emph{region of attraction} (RoA); i.e., the set of states for which the MPC is initially feasible and which are hence stabilized by the MPC.

%%-------------------------------------------------------------------------------

\subsection{Results for DRMPC}

\begin{figure}[h!]
    \centering
    \includegraphics[width=0.4\linewidth]{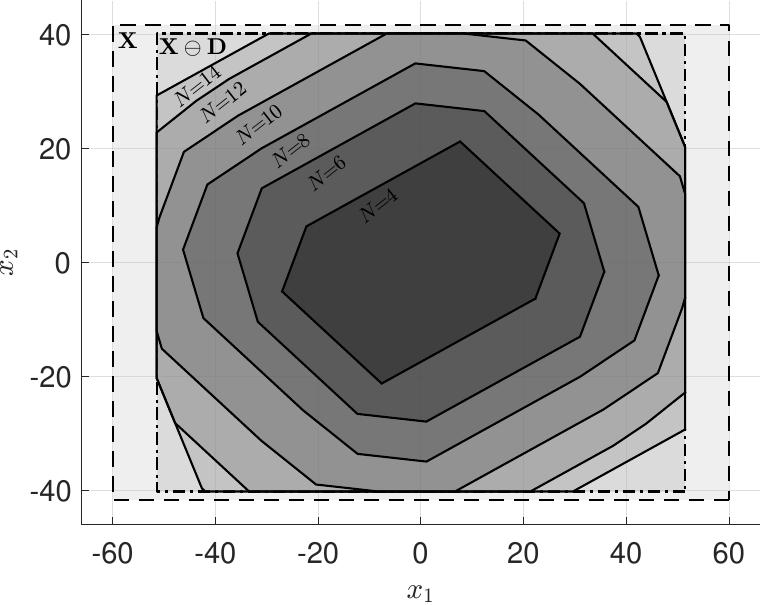}
    \caption{The achived RoAs with the DRMPC for   $N = 4, 6, 8, \dots, 14$. The deadbeat horizon is $M = 3$ and the terminal set is the origin.}
    \label{fig:DoA_DRMPC}
\end{figure}

For the Robust Case, Figure~\ref{fig:DoA_DRMPC} compares the respective RoA 
for different horizon lengths $N = 4, 6, 8, 10, 12, 14$. As expected with the zero terminal constraint, the RoA is small for short horizons.
However, increasing the prediction horizon leads to a reasonably fast growth of the RoA, which mitigates the limitation imposed by the zero terminal constraint. Beyond $N = 14$, further increasing $N$ does not result in a significant improvement of the RoA. 

In fact, the state constraint set is restricted to $\mathbb{X}\ominus \BD$, and is tightened further over the prediction horizon, as defined in~(\ref{Equ:AuxMatrices}a-f). Therefore, the maximum achievable RoA based on the DRMPC must lie within the set $\mathbb{X} \ominus \BD$; see Figure~\ref{fig:DoA_DRMPC}. This is a consequence of the procedure presented in Section~\ref{Sec:AdditiveUncertainty} for converting parametric uncertainty into an additive disturbance in the Robust Case, depending on the size of the resulting disturbance set~$\BD$, which constitutes the main limitation of this procedure. Interestingly, this limitation does not arise in the LPV Case, as illustrated in the next section. 

In comparison with a Tube MPC approach \cite{FlemingEtAl:2015}, which considered terminal RPI sets with $N=6$, it is clear that the RoA achieved with the proposed DRMPC is smaller than the one with a terminal RPI set \cite{FlemingEtAl:2015}. However, slightly increasing $N$ allows the DRMPC to achieve a larger RoA, even with the zero terminal set (though it is still confined within the set $\mathbb{X} \ominus \BD$).

In terms of computational complexity, for $N=6$, the number of decision variables in the Tube MPC \cite{FlemingEtAl:2015} is 85, with 342 inequality constraints and 16 equality constraints. In contrast, the DRMPC formulation requires only 6 decision variables, 32 inequality constraints, and 16 equality constraints. This indicates a significant reduction in computational complexity, even with a longer prediction horizon. For $N=14$, the DRMPC optimization problem consists of 14 decision variables, 82 inequality constraints, and 32 equality constraints.

%%-------------------------------------------------------------------------------
\subsection{Results for LPV-DRMPC}

In the LPV Case, a realization of the nominal future values of $\theta$, i.e., $\bar{\theta}$, should be available according to Assumption~L.\ref{The:Measurement-theta}.Therefore, the resulting RoA depends on the specific nominal values $\bar{\theta}$ considered.
%In the literature on LPV-MPC (e.g., \cite{}), the RoA is commonly presented based on a single realization of $\theta$, and the impact of that realization on the resulting RoA is often overlooked.  This example highlights that aspect.
To evaluate the impact of the realization of $\bar{\theta}$  on the resulting RoA, different realizations of the nominal values under varying conditions are hence considered. The resulting RoAs are compared in terms of their volume for different values of the prediction horizon $N$ and different bounds on $\Delta\theta$, denoted by $\Delta\theta_{\max}$, such that $\Delta\theta \in [-\Delta\theta_{\max}, \Delta\theta_{\max}]$, which defines the set $\Theta_\Delta$, see Fig.~\ref{Fig:schtube}. Finally, the results for the LPV-DRMPC are compared to a state-of-the-art Tube LPV-MPC scheme \cite{Abbas:2019}.
 
 %Therefore, we consider three MPC-LPV controllers based on different ways of realizing $\theta$, and we compare them in terms of the achieved RoA, specifically its volume, for different values of the prediction horizon $N$. Then, we compare our approach with  other approaches in \cite{} and \cite{}. 

%\subsubsection{Design 1}

In the following, we consider random sequences of the nominal values of $\theta$, but with a restricted rate of variation bounded by $\Delta\theta_{\max}$, as is typically the case in practice.
% that does not exceed its uncertainty bound $\Delta\theta$, which defines the associated set $\Theta_\Delta$. 
%This is realistic, as in many applications the trajectory of $\theta$ is subject to a bounded rate of change. 
A total of 144 realizations have been generated, each comprising a sequence of 12 points, and used to construct the RoAs for different values of $N = 4, 5, \dots, 12$ per sequence. Based on the volume of the RoA, we perform the following analysis: First, we compute the ratio of the average  RoA volume to maximum and to minimum RoA volume, denoted $r_{\rm avmax}$ and $r_{\rm avmin}$, respectively, across the 144 realizations for each $N$. Second, we compute the percentage increase in the RoA volume from $N$ to $N+1$, for all $N = 4, 5, \dots, 11$, for each sequence. We then compute the average and standard deviation of this percentage increase across all 144 sequences, denoted as $r_{\rm inc}$. This analysis is repeated       
 for $\Delta\theta_{\max} = 0.1, 0.2, \dots, 0.5$.
%, using the same 144 sequences, taking into account the limits of $\theta$ according to the  parameter set $\Theta$.
All results are presented in Table~\ref{tab:Res-DRMPC-LPV}. 

Now we discuss the results depicted in the table. It turns out that for $N = 10, 11, 12$, no noticeable change in the volume or shape of the RoA is observed across all sequences and considered values of $\Delta\theta_{\max}$, and they are therefore omitted accordingly.
With respect to the values of $r_{\rm avmax}$ and $r_{\rm avmin}$, for $\Delta\theta_{\max} = 0.1$, we observe from Table~\ref{tab:Res-DRMPC-LPV} that at $N = 4, 5, 6$, the average RoA volume is approximately 1.15–1.18 times the minimum RoA volume and about 0.81–0.87 times the maximum RoA volume. This suggests that with small $N \leq 6$, the size of the RoA is somewhat affected by the realization of $\bar{\theta}$, but this effect is not significant. In contrast, for large $N$, i.e., $N \geq 7$, $r_{\rm avmin}$  
approaches $1.02$–$1.09$, and $r_{\rm avmax}$ approaches $0.95$–$0.99$, indicating that the average RoA becomes  
very close to both the minimum and maximum RoA. Hence, for larger prediction horizons, changing the realization of $\bar{\theta}$ has little effect on the size of the RoA.
Now we consider $r_{\rm inc}$. For $\Delta\theta_{\max} = 0.1$, the RoA increases by approximately $79 \pm 2.8\%$ when going from $N=3$ to $N=4$, about $55 \pm 4.7\%$ from $N=4$ to $5$, $21 \pm 3.7\%$ from $N=5$ to $6$, and less than $1\%$ for $N \geq 8$. This trend indicates that most of the improvement in RoA size occurs when increasing $N$ from small values, and the effect of the realization of $\bar{\theta}$ is minor, as reflected by the small standard deviation in these results.

\begin{figure}[ht]
    \centering
        \includegraphics[width=0.3\linewidth]{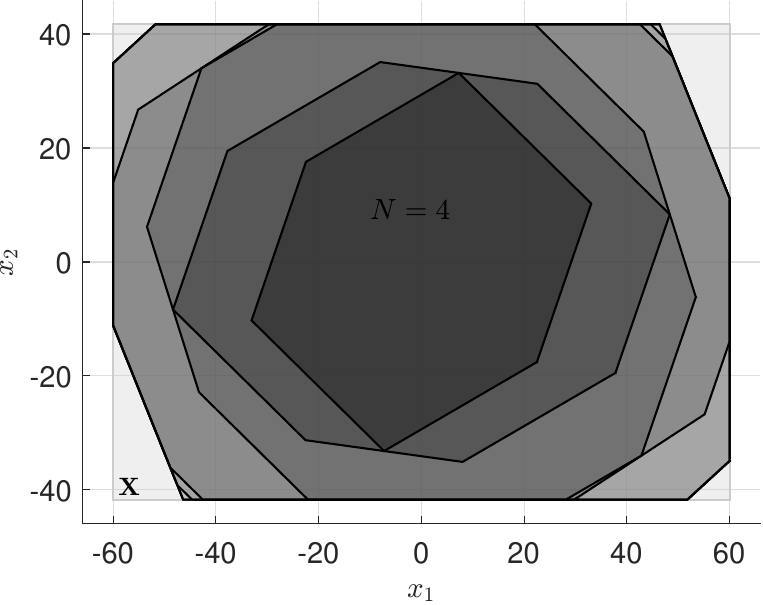}
        \includegraphics[width=0.3\linewidth]{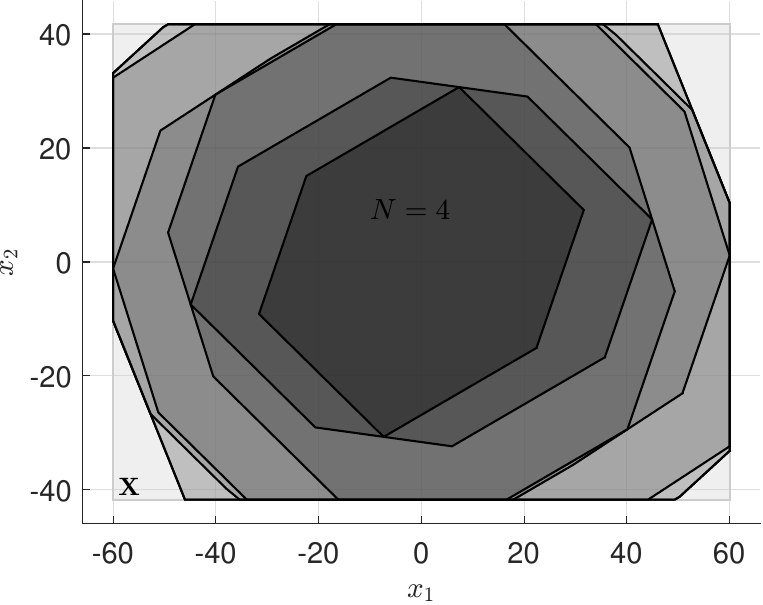}
        \includegraphics[width=0.3\linewidth]{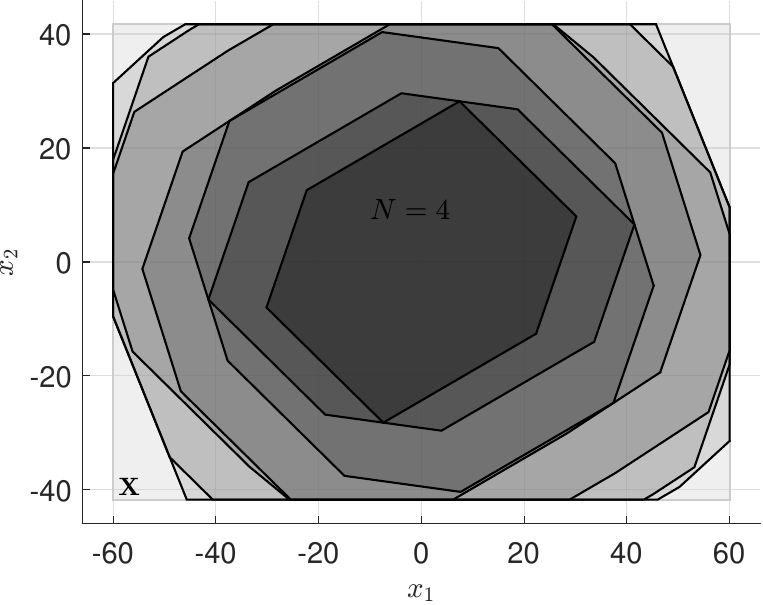}
    \caption{ Representative RoAs corresponding to the sequence shown in Fig~\ref{fig:theta_realization} with different $n = 4, 5, \dots, 10$ and the uncertainty bounds: (left) $\Delta\theta_{\max} = 0.1$, (middle) $\Delta\theta_{\max} = 0.3$, (right) $\Delta\theta_{\max} = 0.5$.}
    \label{fig:roa_lpv}
\end{figure}

A similar pattern is observed across all $\Delta\theta_{\max} = 0.1, 0.2, \dots, 0.5$ for $r_{\rm avmin}$, $r_{\rm avmax}$, and $r_{\rm inc}$. This indicates that the general behavior of the proposed LPV-DRMPC is not significantly affected by the uncertainty bound $\Delta\theta_{\max}$, although increasing $\Delta\theta_{\max}$ generally reduces the size of the RoA. This is illustrated by the representative RoA results shown in Figure~\ref{fig:roa_lpv}, which corresponds to one sequence of $\bar{\theta}$ used to construct Table~\ref{tab:Res-DRMPC-LPV} with  $\Delta\theta_{\max}=0.1,0.3,0.5$. It is observed that the RoAs in Fig.~\ref{fig:roa_lpv}(a) for $\Delta\theta_{\max} = 0.1$ and $N = 8, 9, 10$ are not distinguishable, whereas for $\Delta\theta_{\max} = 0.3$ in Fig.~\ref{fig:roa_lpv}(b), this is the case for $N = 9, 10$.
It is also shown that the largest RoAs across all cases are nearly similar, which indicates that the largest achievable RoA under increased uncertainty bounds  can be maintained by slightly increasing the prediction horizon length.

\begin{figure}[h!]
        \centering
        \scalebox{1}{\input{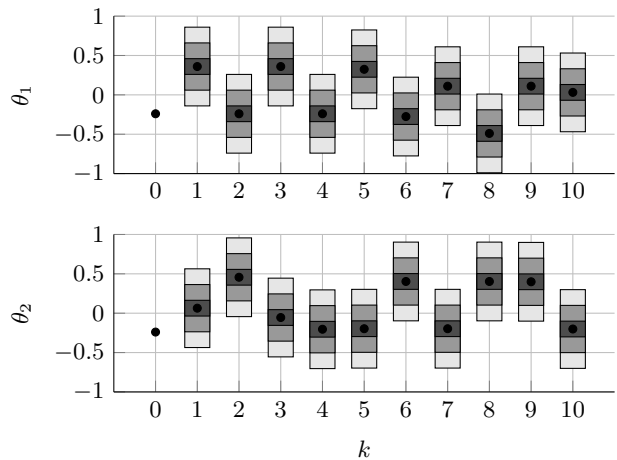}} 
        \caption{The realization of $\theta_1$ and $\theta_2$ with their associated uncertainty bounds $\Delta\theta_{\max}$ used for producing the RoAs.}
        \label{fig:theta_realization}
    \end{figure}

In the following, we compare our approach with the tube-based LPV-MPC  in \cite{Abbas:2019}, which comes with theoretical guarantees on recursive feasibility and stability.
This method does not rely on given future nominal values of $\theta$; instead, it constructs admissible values of $\theta$ from specified bounds on its rate of change and defines the tubes accordingly. These tubes are constructed online and used for online constraint tightening.
Unlike our  approach, it incorporates an MPI set as a terminal constraint based on an LQR design.
For a fair comparison, we consider the tube-based approach solved using given nominal values $\bar{\theta}$ with uncertainty bounds around them, similar to our method. In particular, we use the nominal realization shown in Fig.~\ref{fig:theta_realization}, with different values of $\Delta\theta_{\max}$.

For a quantitative comparison, Table~\ref{tab:comp-with-tube-based} shows the excess volume (in percent) of the RoAs obtained with LPV-DRMPC ---using either a terminal set at the origin ($x_N = 0$) or a terminal set given by an LQR-based RPI set $\mathbb{X}_f$ (i.e., $x_N \in \mathbb{X}_f$)--- over those obtained with the tube-based LPV-MPC, for different uncertainty bounds $|\theta|$ and prediction horizons $N$.

%For a quantitative comparison, Table~\ref{tab:comp-with-tube-based} shows the excess volume of  RoAs of  LPV-DRMPC with a terminal set at the origin, i.e., $x_N = 0$, and with a terminal set as a RPI set $\mathbb{X}_f$ based on an LQR, i.e., $x_N\in \mathbb{X}_f$, over that of the tube-based LPV-MPC and we consider with different uncertainty bounds $\|\theta\|$ and different prediction horizons $N$, where the same realization in Fig.~\ref{fig:theta_realization}s is used for all cases.

It is demonstrated from Table~\ref{tab:comp-with-tube-based} that the performance of the proposed LPV-DRMPC scheme
with $x_N=0$ and a small prediction horizon is worse than that of the tube-based approach, as indicated by the large negative values. However, by slightly increasing $N$ to 9 or 10, the performance of the LPV-DRMPC with $x_N=0$  becomes almost comparable to that of the tube-based LPV-MPC. On the other hand, the  LPV-DRMPC with the terminal set $\mathbb{X}_{\rm f}$ performs nearly similarly to the tube-based LPV-MPC, even for small $N$. It is noticeable that avoiding the complexity of computing a terminal set (by considering the zero terminal condition) requires slightly increasing the prediction horizon
to recover the loss in performance. It is also worth mentioning that the computational complexity of the proposed approach is,
in general, significantly less than that of the tube-based MPC, as it avoids the online construction of tubes and the constraints tightening.

\begin{figure}[h!]
    \centering
        \includegraphics[width=0.3\linewidth]{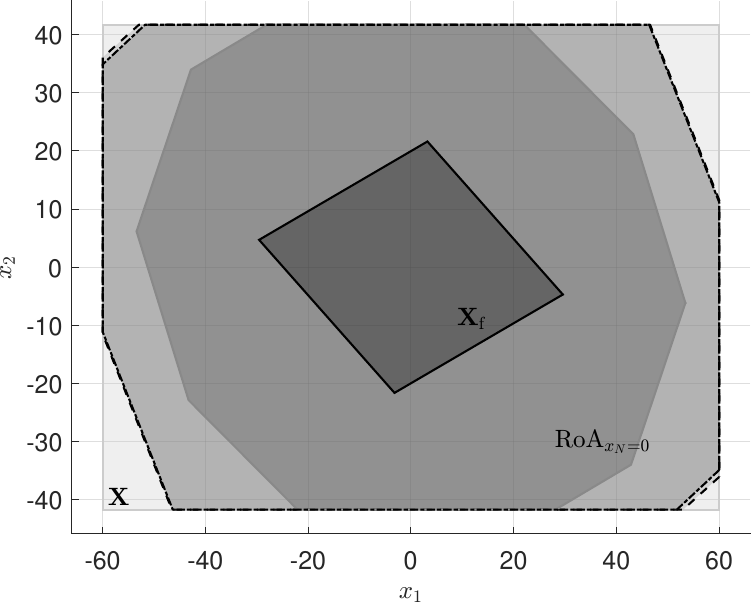}
        \includegraphics[width=0.3\linewidth]{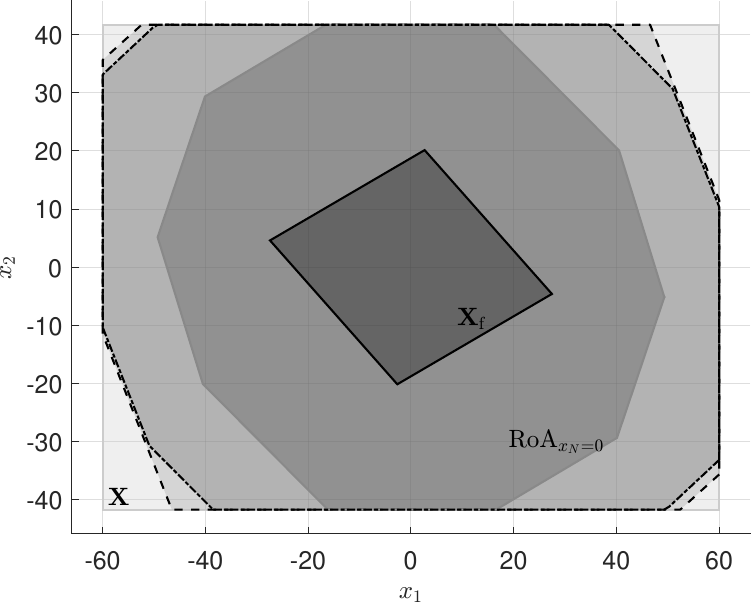}
        \includegraphics[width=0.3\linewidth]{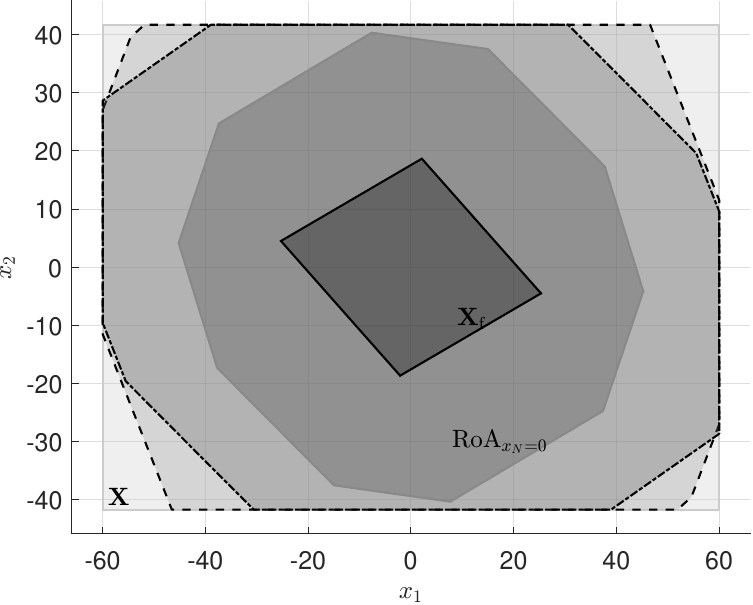}
    \caption{Comparison with tube-based: Representative RoAs corresponding to  the sequence shown in Fig~\ref{fig:theta_realization} with  $N = 6$ and the uncertainty bounds: (left) $\Delta\theta_{\max}= 0.1$, (middle) $\Delta\theta_{\max}= 0.3$, (right) $\Delta\theta_{\max}= 0.5$.}
    \label{fig:comp-with-tube-based-N6}
\end{figure}

\begin{figure}[h!]
    \centering
        \includegraphics[width=0.3\linewidth]{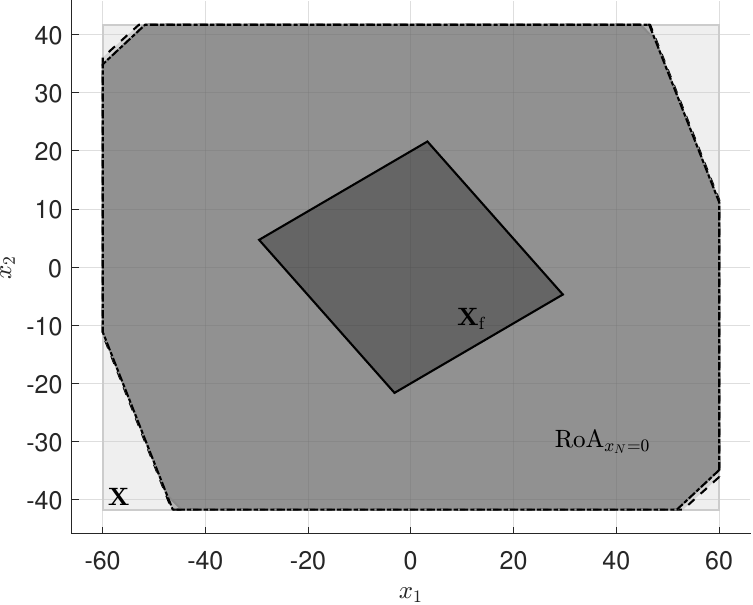}
        \includegraphics[width=0.3\linewidth]{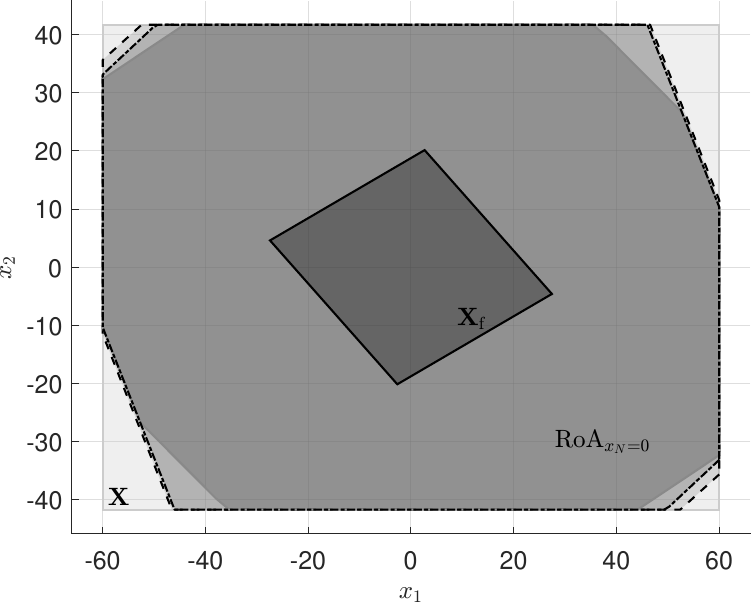}
        \includegraphics[width=0.3\linewidth]{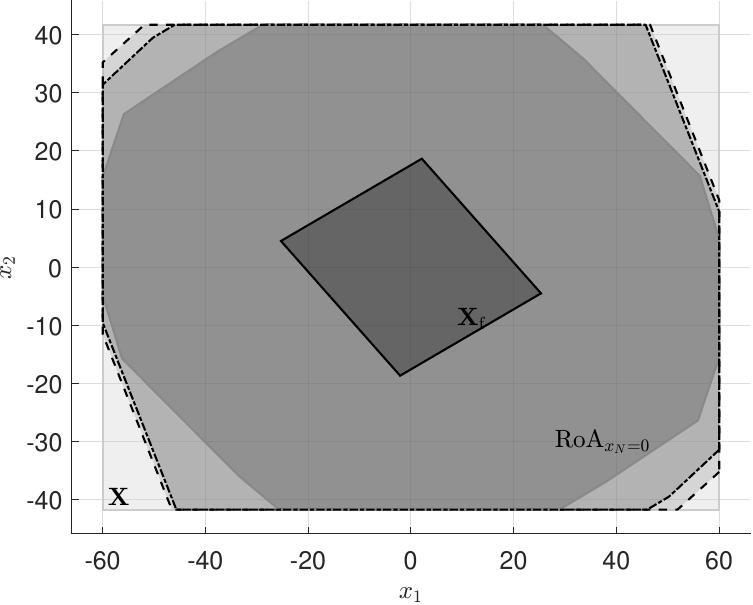}
    \caption{Representative RoAs corresponding to the sequence shown in Fig~\ref{fig:theta_realization} with   $N =8$ and the uncertainty bounds: (left) $\Delta\theta_{\max}= 0.1$, (middle) $\Delta\theta_{\max}= 0.3$, (right) $\Delta\theta_{\max}= 0.5$.}
    \label{fig:comp-with-tube-based-N8}
\end{figure}

As representative results from Table~\ref{tab:comp-with-tube-based}, Figs.~\ref{fig:comp-with-tube-based-N6} and \ref{fig:comp-with-tube-based-N8} depict the RoAs for $N = 6$ and $N = 8$, respectively, along with the  terminal set $\mathbb{X}{\rm f}$, for the LPV-DRMPC results with both the zero terminal set and $\mathbb{X}{\rm f}$, in comparison with the tube-based LPV-MPC and under  different uncertainty  bounds $\Delta\theta_{\max}$.

\begin{table}[htbp]
\centering
\renewcommand{\arraystretch}{1.2}
\setlength{\tabcolsep}{4pt}
\begin{tabular}{|c||c||*{6}{c|}} 
\hline
{$\Delta\theta_{\max}$}   & \multirow{3}{*}{\begin{tabular}{c} \\ \end{tabular}} 
  & \multicolumn{6}{c|}{Horizon length $N$} \\ 
\cline{3-8}
  & & 4 & 5 & 6 & 7 & 8 & 9 \\ 
\hline\hline
0.1  &$r_{\rm avmin}$ &1.15 &1.18 &1.15 &1.09 & 1.03&1.02  \\ \cline{2-8}
  &$r_{\rm avmax}$ &0.83 &0.81 &0.87 &0.95 & 0.99&0.99  \\ \cline{2-8}
  &$r_{\rm inc}$ & 79 $\pm$ 2.8& 54.6 $\pm$ 4.7& 20.8 $\pm$ 3.7& 6.1 $\pm$ 2.3&0.2 $\pm$ 0.3 &     0  \\ \cline{2-8}
  \hline\hline
0.2  &$r_{\rm avmin}$ &1.15 &1.17 &1.16 &1.12 & 1.04&1.02  \\ \cline{2-8}
  &$r_{\rm avmax}$ &0.83 &0.81 &0.86 &0.93 & 0.98&0.99  \\ \cline{2-8}
  &$r_{\rm inc}$ &75.7$\pm$2.3 &  56.5$\pm$3.9 &       24.3$\pm$3.4 &       8.8$\pm$2.8 &      0.6$\pm$0.8  &    0  \\ 
  \hline\hline
0.3  &$r_{\rm avmin}$ &1.14 &1.16 &1.16 &1.13 & 1.06&1.02  \\ \cline{2-8}
  &$r_{\rm avmax}$ &0.83 &0.81 &0.85 &0.91 & 0.97&0.99  \\ \cline{2-8}
  & $r_{\rm inc}$&72.3$\pm$2.0&57.6$\pm$3.1& 27.1$\pm$2.9&  12.5$\pm$3.4&1.9$\pm$1.6& 0 \\ 
  \hline \hline
0.4  &$r_{\rm avmin}$ &1.13 &1.15 &1.16 &1.13 & 1.08&1.03  \\ \cline{2-8}
  &$r_{\rm avmax}$ &0.82 &0.80    &0.83 &0.87& 0.94&0.98  \\ \cline{2-8}
  & $r_{\rm inc}$& 68.8$\pm$1.6 &57.6$\pm$2.1&29$\pm$2.5& 16.9$\pm$3.8&4.6$\pm$2.3&0.2$\pm$0.3  \\ 
  \hline \hline
0.5  &$r_{\rm avmin}$ &1.11 &1.13 &1.14 &1.13 & 1.09&1.05  \\ \cline{2-8}
  &$r_{\rm avmax}$ &0.80 &0.80 &0.81 &0.84 & 0.91&0.97  \\ \cline{2-8}
  & $r_{\rm inc}$& 65.2$\pm$1.3&  56.3$\pm$1.5& 30.7$\pm$ 1.9  &21.4$\pm$ 3.8  & 8.4$\pm$2.6&           1.3$\pm$1.2  \\ 
  \hline
\end{tabular}
\caption{Results of the LPV-DRMPC for different realizations of the nominal parameter $\bar{theta}$}
\label{tab:Res-DRMPC-LPV}
\end{table}

\begin{table}[htbp]
\centering
\renewcommand{\arraystretch}{1.2}
\setlength{\tabcolsep}{4pt}
\begin{tabular}{|c||l||c|c|c|c|c|} 
\hline
{$\Delta\theta_{\max}$} &  & \multicolumn{5}{c|}{Horizon length $N$} \\ 
\cline{3-7}
 & & 6 & 7 & 8 & 9 & 10 \\ 
\hline\hline
0.1  &  $x_{N}=0$ &$-24.05\%$ & $-6.08\%$  & $-0.33\%$  &$-0.28\%$  &  $-0.28\%$ \\ \cline{2-7}
      & $x_{N}\in\mathbb{X}_{\rm f}$   &  $-0.28\%$ & $-0.28\%$ & $-0.28\%$ & $-0.28\%$ &  $-0.28\%$ \\ 
\hline\hline
0.2  & $x_{N}=0$ & $-28.67\%$  &$-9.65\%$  & $-1.07\%$ & $-0.53\%$ & $-0.53\%$  \\ \cline{2-7}
      & $x_{N}\in\mathbb{X}_{\rm f}$   & $-0.74\%$ & $-0.53\%$ & $-0.53\%$ & $-0.53\%$ & $-0.53\%$  \\ 
\hline\hline
0.3  & $x_{N}=0$ & $-33.49\%$ & $-15.34\%$  &$-3.03\%$  & $-0.79\%$ & $-0.79\%$  \\ \cline{2-7}
      & $x_{N}\in\mathbb{X}_{\rm f}$   & $-1.65\%$  &$-0.79\%$  & $-0.79\%$ & $-0.79\%$ &  $-0.79\%$ \\ 
\hline\hline
0.4  & $x_{N}=0$ & $-38.85\%$ & $-21.86\%$ &$-6.50\%$  & $-1.16\%$ & $-1.08\%$  \\ \cline{2-7}
      & $x_{N}\in\mathbb{X}_{\rm f}$   & $-3.35\%$ & $-1.08\%$ & $-1.08\%$ & $-1.08\%$ &  $-1.08\%$ \\ 
\hline\hline
0.5  & $x_{N}=0$ & $-44.58\%$ & $-28.57\%$ & $-11.76\%$ & $-3.11\%$ & $-1.40\%$  \\ \cline{2-7}
      & $x_{N}\in\mathbb{X}_{\rm f}$  & $-5.94\%$ & $-1.56\%$ &$-1.40\%$  & $-1.40\%$ &  $-1.40\%$ \\ 
\hline
\end{tabular}
\caption{Percentage excess volume of the RoA of LPV-DRMPC in both cases, $x_{N|0}=0$ and $x_{N|0}\in\mathbb{X}_{\rm f}$, over tube-based LPV-MPC.}
\label{tab:comp-with-tube-based}
\end{table}

%\begin{figure}[ht]
%    \centering
%%    \includegraphics[width=0.4\linewidth]{DoA_LPVcase.pdf}
 %   \caption{ROA LPV Case case $N=4:8$.}
 %   \label{fig:roa_lpv}
%\end{figure}

\newlength\figH 
\newlength\figW
\setlength\figH{5.cm} 
\setlength\figW{7.cm}

%\backmatter

%%%%%%%%%%%------------------------------------

\section{Conclusion}\label{Sec:Conclusion}

This paper has proposed a DRMPC for linear systems and for a class of LPV systems subject to parametric uncertainties and additive disturbances. For the LPV system, a nominal trajectory of its future varying parameter is assumed to be given with uncertainty bounds, while its current value is known.
To design a computationally efficient MPC scheme, the contributions of the parametric uncertainty have been lumped together with the additive disturbance, and a worst-case bound for this quantity has been obtained accordingly. Thus, the DRMPC is readily applicable and allowed completely extinguishing the additive disturbances and the effects of uncertainties within the deadbeat horizon.
To reduce online computational complexity with a slight increase in conservatism, the deadbeat inputs are precalculated during the offline design phase, allowing all constraint tightening to be performed offline, leading to a more practical approach for applications in practice.
One of the main advantages of LPV-DRMPC, however, is that no calculation of any PI or RPI sets or terminal costs is needed, as the terminal constraint set is taken as the origin. In comparison with existing approaches, where a tube-based LPV-MPC scheme has been considered, the performance of LPV-DRMPC turned out to be competitive, and the loss of performance could be almost fully recovered by slightly extending the prediction horizon.
This implies that it is straightforward to extend the LPV-DRMPC approach to nonlinear systems embedded in the LPV framework. Moreover, the extension of the entire approach is plausibly worthwhile for data-driven predictive control based solely on system data, offering theoretical guarantees without the need for terminal conditions. The details of these ideas will be explored in future research.

%%%%%%%%%%%------------------------------------

\section*{Acknowledgments}

Research leading to these results has received funding from the German Research Foundation (DFG) under grants no.\ 460891204 and no.\ 419290163.

%\bmsection*{Conflict of interest}

%The authors declare no potential conflict of interests.
%%%%%%%%%%%------------------------------------
\bibliographystyle{h-physrev3}
\bibliography{bibcontr,bibeng,bibmath}

\end{document}